  \documentclass[12pt]{amsart}
  \usepackage{latexsym}
  \usepackage[all]{xy}
  \usepackage{amsfonts}
  \usepackage{amsthm}
  \usepackage{amsmath}
  \usepackage{amssymb}
  \usepackage{hyperref}

 \numberwithin{equation}{section}
 \usepackage{mathptmx}
 \usepackage{color}

  \def\sw#1{{\sb{(#1)}}}
  
  \def\su#1{{\sp{[#1]}}}
  \def\suc#1{{\sp{(#1)}}}

  \def\<{{\langle}}
  \def\>{{\rangle}}

  \def\note#1{{}}

  \def\note#1{}
  \def\M{{\bf M}}

  \def\lrhom#1#2#3#4{{{\rm Hom}\sb{#1, #2}(#3,#4)}}
   \def\lRhom#1#2#3#4{{{\rm Hom}\sb{#1-}\sp{-#2}(#3,#4)}}

  \def\rend#1#2{{{\rm End}\sb{-#1}(#2)}}

  \def\lrend#1#2#3{{{\rm End}\sb{#1,#2}(#3)}}

  \def\Bihom#1#2#3#4{{{\rm Hom}\sp{#1,#2}(#3,#4)}}

  \def\beq{\begin{equation}}
  \def\eeq{\end{equation}}

  \def\ot{{\otimes}}
  
  \def\Hom{\mbox{\rm Hom}\,}

  \newcommand{\Ra}{\Rightarrow}

  \def\Aro{{{}^A\!\varrho}}
  \def\roA{{\varrho^A}}

  \newcommand{\cC}{\mathcal{C}}
  \def\M{{\bf M}}
  \def\ot{{\otimes}}
  \def\hH{\mathcal H}
  \def\Aro{{{}^A\!\varrho}}
  \def\roA{{\varrho^A}}
  \def\ui{^{(1)}}
  \def\uii{^{(2)}}
  \def\di{_{(1)}}
  \def\dii{_{(2)}}
  \def\tz{^{[0]}}
  \def\ti{^{[1]}}
  \def\bz{_{[0]}}
  \def\bi{_{[1]}}

  \def\cp{\raise1pt\hbox{${\scriptstyle{\#}}$}}
\def\stac#1{\raise-.2cm\hbox{$\stackrel{\displaystyle\otimes}{\scriptscriptstyle{#1}}$}}
\def\sstac#1{\otimes_{#1}}
\def\cten#1{\raise-.2cm\hbox{$\stackrel{\displaystyle\widehat{\otimes}}
{\scriptscriptstyle{#1}}$}}

\def\conv#1{#1^\textsf{c}}
\def\lconv#1{#1^\textsf{w}}
\def\cL{\mathcal L}
\def\cR{\mathcal R}
\def\Conv#1#2{{\textsf{Conv}}(#1,#2)}
  \newcounter{zlist}
  \newenvironment{zlist}{\begin{list}{(\arabic{zlist})}{
  \usecounter{zlist}\leftmargin2.5em\labelwidth2em\labelsep0.5em
  \topsep0.6ex
  \parsep0.3ex plus0.2ex minus0.1ex}}{\end{list}}

  \newcounter{blist}
  \newenvironment{blist}{\begin{list}{(\alph{blist})}{
  \usecounter{blist}\leftmargin2.5em\labelwidth2em\labelsep0.5em
  \topsep0.6ex 
  \parsep0.3ex plus0.2ex minus0.1ex}}{\end{list}}
\def\stac#1{\raise-.2cm\hbox{$\stackrel{\displaystyle\otimes}{\scriptscriptstyle{#1}}$}}
\def\sstac#1{\otimes_{#1}}

  \newtheorem{proposition}{Proposition}[section]
  \newtheorem{lemma}[proposition]{Lemma}
  
  \newtheorem{corollary}[proposition]{Corollary}
  \newtheorem{theorem}[proposition]{Theorem}

  \theoremstyle{definition}
  \newtheorem{definition}[proposition]{Definition}
  \newtheorem{example}[proposition]{Example}

  \theoremstyle{remark}
  \newtheorem{remark}[proposition]{Remark}

  \newcounter{c}
  
  \newcommand{\etyk}[1]{\vspace{-7.4mm}$$\begin{equation}\Label{#1}
  \addtocounter{c}{1}}
  \renewcommand{\]}{\ifnum \value{c}=1 $$\else \end{equation}\fi}
\begin{document}

 \title{Cleft extensions of Hopf algebroids}
 \author{Gabriella B\"ohm}
 \address{Research Institute for Particle and Nuclear Physics, Budapest,
 \newline\indent H-1525
 Budapest 114, P.O.B.\ 49, Hungary}
  \email{G.Bohm@rmki.kfki.hu}
  \author{Tomasz Brzezi\'nski}
  \address{ Department of Mathematics, University of Wales Swansea,
  Singleton Park, \newline\indent  Swansea SA2 8PP, U.K.}
 \email{T.Brzezinski@swansea.ac.uk}
    \date{October 2005}
  \begin{abstract}
 The notions of a cleft extension and a cross product with
  a Hopf algebroid are introduced and studied. In particular it is
  shown that an extension (with a Hopf algebroid $\hH= (\hH_L,\hH_R)$) is
  cleft if and
  only if it is $\hH_R$-Galois and has a normal basis property relative
  to the base
  ring $L$ of $\hH_L$. Cleft extensions are identified
  as  crossed products with invertible cocycles. The relationship between
   the equivalence
  classes of crossed products and gauge transformations is established.
  Strong connections in cleft extensions are classified and sufficient
  conditions are
  derived for the Chern-Galois characters to be independent on the choice
  of strong connections. The results concerning cleft extensions and crossed
  product are then extended to the case of {\em weak} cleft extensions of Hopf
  algebroids hereby defined.
   \end{abstract}
  \maketitle

\section{Introduction}
Cleft extensions of algebras by a Hopf algebra, or cleft Hopf comodule
algebras, are one of the simplest
and best known examples of Hopf-Galois extensions. Indeed, by
\cite[Theorem~9]{DoiTak:cle} a Hopf-Galois extension with the normal basis
property is necessarily a cleft extension. With geometric interpretation of
Hopf-Galois extensions over fields 
as non-commutative principal bundles, cleft extensions
can be understood as such principal bundles that
every associated bundle is trivial. Motivated by examples coming from
non-commutative differential geometry, the notion of a Hopf-Galois extension
was generalised to a coalgebra-Galois extension
in \cite{BrzMaj:coa}, \cite{BrzHaj:coa}. Subsequently, the notion of a {\em
  cleft
coalgebra extension} was introduced in \cite[p.\ 293]{Brz:mod},
and most comprehensively studied in terms of {\em cleft entwining
structures} in \cite{Abu:Mor}, \cite{CaeVer:Mor}. The latter were
extended further to weak entwining structures in
\cite{AlAl:wcleft}, \cite{AlAl:wcleftgal}. 

The aim of the present paper is to extend the theory of cleft extensions in a
different direction, in the first instance motivated by recent developments
in the theory
of depth 2 and Frobenius ring extensions \cite{Kad:D2gal}, \cite{Kad:norm},
\cite{KadSzl:D2bgd}, in long term motivated by the increasing interest in
Galois-type extensions with Hopf algebroid symmetries \cite{Bohm:gal},
\cite{BalSzl:fgal}. Thus we introduce
and study basic properties of Hopf algebroid cleft extensions. Very much as
cleft extensions for Hopf algebras are an example and a testing ground for
more general Hopf-Galois extensions, also  Hopf algebroid cleft extensions
provide one with a useful tool (or a toy model) for more general Hopf algebroid
extensions. In particular, as announced in \cite{BohBrz:str}, cleft
extensions for Hopf algebroids give a concrete illustration to the
relative Chern-Galois theory. In fact the current paper can be considered as
a sequel to \cite{BohBrz:str} in which the ideas and results, announced in
a few examples, are developed in detail and further extended. Specifically, in
 Section~\ref{sec.chern}, sufficient conditions for the existence of
  (strong connection independent) relative Chern-Galois characters
in Hopf algebroid cleft extensions are stated.

The construction of Hopf algebroid cleft extensions, although
motivated by similar ideas, is significantly different from that of
cleft Hopf algebra (or coalgebra) extensions. One should remember
that a Hopf algebroid involves two different coring (and bialgebroid)
structures on the same
$k$-module. The
interplay between these intricate structures is an immanent feature of
Hopf algebroid extensions. This is already present in the notion of
 a convolution inverse (cf.\ Definition~\ref{def.conv.inv}), which relates
 two coring structures on the same $k$-module, but is perhaps most
 significant in the characterisation of cleft extensions in terms of
 the Galois and normal basis properties (cf.\ Theorem~\ref{thm.main.cleft}):
 a cleft $\hH$-extension is a Galois extension with
 respect to the {\em right} bialgebroid $\hH_R$
 but it has a normal basis property with respect
 to the {\em left} bialgebroid $\hH_L$.

 In the standard Hopf algebra theory, cleft extensions of Hopf algebras
 are examples of {\em crossed products} with Hopf algebras: indeed a cleft
 extension is the same as a crossed product with an invertible cocycle
 (cf.\ \cite[Theorem~11]{DoiTak:cle},
\cite[Theorem~1.18]{BlaMon:cro}). Motivated by this correspondence,
we also develop a general theory of crossed products with bialgebroids and
Hopf algebroids. In particular this involves developing the notions of
a {\em measuring} and a {\em 2-cocycle}, while to relate different crossed
products one needs to give meaning to
{\em gauge transformations} and {\em equivalent crossed products}. In parallel
 to
the bialgebra case, we show in Theorem~\ref{thm.equiv}
that two crossed products are equivalent if and only
if one is a gauge transform of the other. We then identify cleft extensions
of Hopf algebroids with crossed products with invertible cocycles (cf.\
Theorem~\ref{thm.cross.cleft} and Theorem~\ref{thm:crossed}). A generalisation
 of
this theory to the case of {\em weak crossed products} is then outlined in
Appendix.

Finally, we would like to indicate that the cleft extensions of the present
paper can be placed in a broader context. 
A (weak) entwining structure $(A,D,\psi)$ determines a 
coring extension $D$ of the canonical $A$-coring $\cC_\psi$ 
{\color{blue}
that is {\em pure} in the sense that it satisfies assumptions of
\cite[Theorem~2.6~(2), (arXiv version)]{Brz:not}; see \cite[Corrigendum,
  Definition~1]{BohmVer:cleftex}}.  
The cleft property of an entwining structure can
be formulated as a feature of $A$ as an entwined module (i.e.\ a
$\cC_\psi$-comodule). Although it is not possible to find a cleft
entwining structure behind a cleft extension $B\subseteq A$ of a
Hopf algebroid $\hH$ (with left $L$-bialgebroid $\hH_L$ and right
$R$-bialgebroid $\hH_R$), there is still an associated coring
extension. Namely, the constituent $L$-coring in $\hH$ is a right
extension of the $A$-coring $\cC :=A\sstac{R} \hH_R$, such that
$A$ is a $\cC$-comodule. Inspired by this observation,  a unified
approach to all known notions of cleft
extensions in terms of 
{\color{blue}
pure}
coring extensions is  developed in
\cite{BohmVer:cleftex}.

{\color{blue}
After publication of this paper, it turned out that the proofs of 
journal versions of \cite[Theorem 2.6]{Brz:not} and \cite[Proposition
  3.1]{Bohm:gal} contain some unjustified steps. The previous
(published) version of this work heavily relied on the journal version of
\cite[Theorem 2.6]{Brz:not} to which we do not resort in this revision.   
Using informality of the arXiv, for the convenience of the
readers of the earlier version, we write the corrections in blue.}

\medskip

\noindent {\bf Notation.} Throughout this paper we work over an associative
unital
commutative ring $k$. An algebra means an associative unital $k$-algebra. Unit
elements are denoted by $1$ and multiplications by $\mu$ (or by $1_R$, $\mu_R$
if the algebra $R$ needs to be specified). Categories of left,
right, and bimodules for an algebra $R$ are denoted by ${}_R\M$, $\M_R$ and
${}_R\M_R$, respectively. Their hom-sets are denoted by ${\rm Hom}_{R,-}(-,-)$
${\rm Hom}_{-,R}(-,-)$ and ${\rm Hom}_{R,R}(-,-)$, respectively.

Categories of left and right comodules for a coring $\cC$ are denoted by
${}^\cC\M$ and $\M^{\cC}$, respectively. For their hom-sets we write
${\rm Hom}^{\cC,-}(-,-)$ and ${\rm Hom}^{-,\cC}(-,-)$, respectively.
\pagebreak

\section{Preliminaries}\label{sec.prelim}
\subsection{Bialgebroids}\label{sec:bgd}
A bialgebroid \cite{Tak:crR},\cite{Lu:hgd} can be considered as a
generalisation of the notion of a bialgebra to arbitrary
(non-commutative) base algebras. A (left) bialgebroid over a base
algebra $L$ consists of an $L\sstac{k} L^{op}$-ring structure
$(H,\mu,\eta)$ and an $L$-coring structure $(H,\gamma,\pi)$ on the
same $k$-module $H$. Denoting the restriction of the unit map
$\eta:L\sstac{k} L^{op} \to H$ to $L\sstac{k} 1_H$ (the so called
{\em source map}) by $s$ and its restriction to $1_H\sstac{k}
L^{op}$ (the {\em target map}) by $t$, the bimodule structure of
the $L$-coring is given by
$$
lhl'=s(l)t(l')h,\qquad \textrm{for all }l,l'\in L,\ h\in H.
$$
The range of the coproduct is required to be in the {\em Takeuchi product}
$$ H\times_L H :=
\{\ \sum_i h_i\stac{L} h'_i\in H\stac{L} H\ |\
\sum_i h_i t(l)\stac{L} h'_i =\sum_i h_i\stac{L} h'_i s(l)\
  \forall l\in L\ \},
$$
which is, indeed, an algebra by factorwise multiplication. The
following compatibility conditions are required between the
$L\sstac{k} L^{op}$-ring and the $L$-coring structures.
\begin{eqnarray}
\gamma(1_H)&=& 1_H\stac{L} 1_H, \\
\gamma(h h')&=&\gamma(h) \gamma(h'), \label{ax:cpmulti} \\
\pi(1_H) &=& 1_L, \\
\pi\left(h s(\pi(h'))\right)&=&\pi(h h'),\label{ax:augs}\\
\pi\left(h t(\pi(h'))\right)&=&\pi(h h'),\label{ax:augt}
\end{eqnarray}
for all $l\in L$ and $h,h'\in H$.

The $L$-$L$ bimodule structure of the coring, underlying a left bialgebroid,
is defined in terms of the multiplication by $s(l)$ and $t(l)$ on the left.
Right bialgebroids are defined analogously in terms of multiplications on the
right. For more details we refer to \cite{KadSzl:D2bgd}.

Thus a bialgebroid is given by the following data:  $k$-algebras $H$ and
$L$, and  maps $s$ (the source), $t$ (the target), $\gamma$ (the
coproduct)
and $\pi$ (the counit). We write $\cL = (H,L,s,t,\gamma,\pi)$.

Note that if $\cL = (H,L,s,t,\gamma,\pi)$ is a left bialgebroid
then so is the co-opposite $\cL_{cop} =
(H,L^{op},t,s,\gamma^{op},\pi)$, where $L^{op}$ denotes the
algebra that is isomorphic to $L$ as a $k$-module, with
multiplication opposite to the one in $L$, and $\gamma^{op}:H\to
H\sstac{L^{op}} H$, $h\mapsto h\dii \sstac{L^{op}} h\di$ is the
coproduct, opposite to $\gamma:H\to H\sstac{L} H$, $h\mapsto
h\di\sstac{L} h\dii$. The opposite, $\cL^{op}=(H^{op},L,$
$t,s,\gamma,\pi)$ is a right bialgebroid.

\subsection{Hopf algebroids}\label{sec:hgd}
Hopf algebroids with bijective antipodes have been introduced in
\cite{BohmSzl:hgdax}.  In \cite{Bohm:hgdint} the definition was extended by
relaxing the bijectivity of the antipode.

A Hopf algebroid consists of two (a left and a right) bialgebroid structures
on the same total algebra. The source and target maps of the left bialgebroid
$\hH_L=(H,L,s_L,t_L,\gamma_L,\pi_L)$ and of the right bialgebroid
$\hH_R=(H,R,s_R,t_R,\gamma_R,\pi_R)$ are related by the following axioms.
\begin{eqnarray}
&& s_L\circ \pi_L\circ t_R=t_R,\qquad t_L\circ \pi_L\circ s_R=s_R
\quad \textrm{and}\nonumber\\
&& s_R\circ \pi_R\circ t_L=t_L,\qquad t_R\circ \pi_R\circ s_L=s_L.
\label{ax:st}
\end{eqnarray}
These conditions imply that the left coproduct $\gamma_L$ is $R$-$R$ bilinear
and the right coproduct $\gamma_R$ is $L$-$L$ bilinear. Each coproduct is
required to be left and right colinear with respect to the other bialgebroid
structure, i.e. the following axioms are imposed:
\begin{equation}\label{ax:cpcolin}
(\gamma_L \stac{R} H)\circ \gamma_R=(H\stac{L} \gamma_R)\circ \gamma_L
\qquad \textrm{and}\qquad
(\gamma_R \stac{L} H)\circ \gamma_L=(H\stac{R} \gamma_L)\circ \gamma_R.
\end{equation}
An $R\sstac{k}L$-$R\sstac{k}L$ bilinear map $S:H\to H$, i.e.\ a
$k$-linear map, such that
\begin{equation}\label{ax:Sbilin}
S\big(t_L(l) h t_R(r)\big)=s_R(r) S(h) s_L(l)
\quad\textrm{and}\quad
S\big(t_R(r) h t_L(l)\big)=s_L(l) S(h) s_R(r),
\end{equation}
for $r\in R$, $l\in L$ and $h\in H$, is called an {\em antipode} if
\begin{equation}\label{ax:antip}
\mu_H\circ (S\stac{L} H)\circ \gamma_L=s_R\circ \pi_R
\qquad \textrm{and}\qquad
\mu_H\circ (H\stac{R} S)\circ \gamma_R=s_L\circ \pi_L.
\end{equation}
For a Hopf algebroid we use the notation $\hH=(\hH_L,\hH_R,S)$.

Since in a Hopf algebroid there are two coring structures present,
we use two versions of Sweedler's index notation for  coproducts.
For any $h\in H$, we write $\gamma_R(h)=h\ui\sstac{R} h\uii$ (with
upper indices) for the right coproduct and $\gamma_L(h)=h\di
\sstac{L} h\dii$ (with lower indices) for the left coproduct.

\begin{remark}\label{rem:axred}
In the formulation of Hopf algebroid axioms, given in \eqref{ax:st},
\eqref{ax:cpcolin}, \eqref{ax:Sbilin} and \eqref{ax:antip}, the left
bialgebroid $\hH_L$ and the right bialgebroid $\hH_R$ play symmetric 
roles. It turns out, however, that this set of axioms can slightly be
reduced. Namely, the second equality in \eqref{ax:Sbilin} can be derived 
from the other axioms.
This can be seen by the following computation (and its symmetric version,
in which the order of multiplication and roles of $\hH_L$ and $\hH_R$ are
interchanged). For $l\in L$ and $h\in H$,
\begin{eqnarray*}
S\big(h t_L(l)\big)
&=& S\big(t_L(\pi_L(h\dii))h\di t_L(l)\big)
= S\big(h\di t_L(l)\big) s_L(\pi_L(h\dii))\\
&=& S(h\di) s_L\big(\pi_L(h\dii s_L(l))\big)
= S(h\di) {h\dii}\ui s_L(l) S({h\dii}\uii)\\
&=& S({h\ui}\di) {h\ui}\dii s_L(l) S({h\uii})
= s_R(\pi_R(h\ui))s_L(l) S({h\uii})\\
&=& s_L(l)s_R(\pi_R(h\ui))S({h\uii})
= s_L(l) S\big( h\uii t_R(\pi_R(h\ui))\big)\\
&=&s_L(l) S(h).
\end{eqnarray*}
The first equality follows by the fact that $\pi_L$ is counit of $\gamma_L$
and the last one follows since $\pi_R$ is counit of $\gamma_R$. The second and
the penultimate equalities follow by the first equality in
\eqref{ax:Sbilin}. The third equality follows by 
the bialgebroid axiom, requiring that the range of $\gamma_L$ is in the
Takeuchi product $H\times_L H$. The fourth equality follows by \eqref{ax:antip}
and \eqref{ax:st}, as the latter implies -- together with the left
$R$-linearity of 
$\gamma_R$ -- that $\gamma_L(hs_L(l))=h\di s_L(l)\sstac{L} h\dii$. The fifth
equality is a consequence of the right $\hH_R$-colinearity of $\gamma_L$,
i.e. \eqref{ax:cpcolin}. The sixth equality follows by \eqref{ax:antip}, and
the seventh one follows by \eqref{ax:st}, implying $s_L(l)s_R(r)=s_R(r)
s_L(l)$, for $r\in R$ and $l\in L$.
\end{remark}

It is proven in \cite[Proposition 2.3]{Bohm:hgdint} that the
antipode of a Hopf algebroid is both an anti-multiplicative map,
i.e.\ $S(hh')=S(h')S(h)$, for $h,h'\in H$, and an
anti-comultiplicative map, i.e.\ $S(h)\di\sstac{L} S(h)\dii$
$=S(h\uii)\sstac{L} S(h\ui)$ and $S(h)\ui\sstac{R}
S(h)\uii=S(h\dii)\sstac{R} S(h\di)$, for $h\in H$, (note the
appearance of left and right coproducts in both formulae). The
maps
\begin{equation}\label{eq:uv}
\pi_R\circ t_L:L^{op}\to R
\qquad \textrm{and}\qquad
\pi_L\circ s_R:R \to L^{op}
\end{equation}
are inverse algebra isomorphisms.

Note that for a Hopf algebroid $\hH=(\hH_L,\hH_R,S)$ also $\hH^{op}_{cop}=
((\hH_R)^{op}_{cop},(\hH_L)^{op}_{cop},S)$ is a Hopf algebroid. If the
antipode $S$ is bijective then also $\hH_{cop}=
((\hH_L)_{cop},(\hH_R)_{cop},S^{-1})$ and
$\hH^{op}=((\hH_R)^{op},(\hH_L)^{op},S^{-1})$ are Hopf algebroids.

{\bf Convention.} Throughout, whenever it is said `Hopf algebroid $\hH$', it is
meant a Hopf algebroid with all the structure modules and maps as in this
section.

\subsection{Comodule algebras for bialgebroids}\label{sec:combgd}

Let $\cR=(H,R,s,t,\gamma,\pi)$ be a right bialgebroid and $M$ a
{\em right $\cR$-comodule}, that is, a right comodule for the
$R$-coring $(H,\gamma,\pi)$. This means \cite[18.1]{BrzWis:cor}
that $M$ is a right $R$-module and there exists a right $R$-linear
coassociative and counital coaction, $\varrho:M\to M\sstac{R} H$,
$m\mapsto m\tz\sstac{R} m\ti$ (note the upper Sweedler indices
indicating the involvement of a right bialgebroid). By the power
of the bialgebroid structure, $M$ can be equipped with a unique
left $R$-action such that $M$ becomes an $R$-$R$ bimodule and
the range of $\varrho$ is in the Takeuchi product
\begin{equation}\label{eq:MTak}
M\times_R H\colon =
\{\ \sum_i m_i\stac{R} h_i\in M\stac{R} H\ |\
\sum_i rm_i\stac{R} h_i= \sum_i m_i\stac{R} t(r)h_i\
\forall\ r\in R\ \}.
\end{equation}
The left $R$-multiplication in $M$ takes the form
\begin{equation}\label{eq:RacM}
rm=m\tz \pi\big(t(r)m\ti\big)\equiv m\tz \pi\big(s(r)m\ti\big),
\qquad \textrm{for all }r\in R,\ m\in M.
\end{equation}
One checks that any $\cR$-colinear map is $R$-$R$ bilinear. In
particular the coaction satisfies
\begin{equation}\label{eq:robilin}
\varrho(rmr')=m\tz\stac{R} s(r)m\ti s(r'), \qquad \textrm{for all
}r,r'\in R,\ m\in M.
\end{equation}
The category of right $\cR$-comodules is a monoidal category with
a strict monoidal functor to the category ${}_R\M_R$ of $R$-$R$
bimodules \cite[Proposition 5.6]{Schau:strthm}. The $R$-action and
$\cR$-coaction on the tensor product of two comodules $M$ and $N$
are given by
\begin{equation}\label{eq:MNcR}
(m\stac{R}n)\cdot r = m\stac{R} nr,\quad (m\stac{R} n)\tz \stac{R} (m\stac{R}
  n)\ti =
(m\tz\stac {R} n\tz)\stac{R} m\ti n\ti,
\end{equation}
for all $r\in R$, $m\sstac{R} n\in M\sstac{R} N$.

A right $\cR$-{\em comodule algebra} is a monoid in the monoidal category of
right $\cR$-comodules; hence, in particular, it is an $R$-ring.

The $R$-coring $(H,\gamma,\pi)$, underlying a right $R$-bialgebroid $\cR$,
possesses a grouplike element $1_H$. The {\em coinvariants} of a right
$\cR$-comodule $M$ with respect to the grouplike element $1_H$ are the
elements of
$$
M^{co\cR}=\{\ m\in M\ |\ m\tz\stac{R} m\ti=m\stac{R} 1_H\ \}.
$$
It is straightforward to check that if $A$ is a right $\cR$-comodule algebra,
then its coinvariants form a subalgebra $B\colon =A^{co\cR}$. In this case
the algebra extension $B\subseteq A$ is termed a right {\em $\cR$-extension}.

A right $\cR$-extension $B\subseteq A$ is a right {\em $\cR$-Galois} extension
if the canonical map
\begin{equation}\label{eq:Rcan}
{\rm can}_R: A\stac{B} A\to A\stac{R} H, \qquad
a\stac{B} a'\mapsto aa^{\prime [0]}\stac{R} a^{\prime [1]} ,
\end{equation}
is bijective, i.e.\ the $A$-coring $A\sstac{R} H$ with the
coproduct $A\sstac{R}\gamma$, the counit $A\sstac{R}\pi$ and with
the $A$-$A$ bimodule structure $a_1(a\sstac{R}
h)a_2=a_1a{a_2}\tz\sstac{R} h{a_2}\ti$, is a Galois coring
\cite{Kad:D2gal}.

Analogously, a right comodule $N$, with coaction $n\mapsto
n\bz\sstac{L} n\bi$ (note the lower Sweedler indices indicating
the involvement of a left bialgebroid), for the $L$-coring
$(H,\gamma,\pi)$, underlying a left bialgebroid
$\cL=(H,L,s,t,\gamma,\pi)$, can be equipped with a left $L$-action
\begin{equation}\label{eq:LacN}
ln=n\bz\pi(n\bi s(l))\equiv n\bz\pi(n\bi t(l)),
\qquad \textrm{for  all } l\in L, n\in N.
\end{equation}
The category of right $\cL$-comodules is a monoidal category with monoidal
product, the module tensor product over $L^{op}$. The right $L$-action and
$\cL$-coaction on the tensor
product of two $\cL$-comodules $M$ and $N$ are
\begin{equation}\label{eq:L-mon}
(m\stac{L^{op}} n)\cdot l = ml\stac{L^{op}} n, \qquad (m\stac{L^{op}}
n)\bz\stac{L} (m\stac{L^{op}} n)\bi=
(m\bz\stac{L^{op}} n\bz)\stac{L} m\bi n\bi,
\end{equation}
for all $l\in L$,  $m\sstac{L^{op}} n\in M\sstac{L^{op}}
N$. Right 
$\cL$-{\em comodule algebras} are defined as monoids in the
monoidal category of right $\cL$-comodules -- hence they are, in
particular, $L^{op}$-rings. {\em Coinvariants} are defined with
respect to the grouplike element $1_H$. An algebra extension
$B\subseteq A$ is called a right $\cL$-{\em extension}  if $A$ is
a right $\cL$-comodule algebra and $B=A^{co\cL}$. A right
$\cL$-extension $B\subseteq A$ is said to be right $\cL$-{\em
Galois} if the canonical map
\begin{equation}
{\rm can}_L: A\stac{B} A\to A\stac{L} H,\qquad
a\stac{B} a'\mapsto a\bz a'\stac{L} a\bi ,
\end{equation}
is bijective.

Right comodules for a left bialgebroid $\cL$ are canonically identified with
left comodules for the co-opposite bialgebroid $\cL_{cop}$, thus resulting in
a monoidal equivalence ${}^{\cL_{cop}}\M\simeq \M^\cL$. This identification
leads to analogous notions of left comodule algebras, left $\cL$-extensions
and left $\cL$-Galois extensions.

\subsection{Comodule algebras for Hopf algebroids}\label{sec:comhgd}

{\color{blue}
In a Hopf algebroid $\hH$ there are two bialgebroid, hence two
  coring structures 
present. The definition of an $\hH$-comodule in \cite[Definition
3.2]{Bohm:gal} and \cite[Section 2.2]{BalSzl:fgal} refers to them in a
symmetrical way:

\begin{definition}\label{def:hgd_comod}
A {\em right comodule} of a Hopf algebroid $\hH$ is a right 
$L$-module as well as a right $R$-module $M$, together with a right
coaction $\varrho_R: M \to M\ot_R H$ of the constituent right bialgebroid
$\hH_R$ and a right coaction $\varrho_L:M\to M\ot_L H$ of the constituent
left bialgebroid $\hH_L$, 
such that $\varrho_R$ is an $\hH_L$-comodule map and $\varrho_L$ is an
$\hH_R$-comodule map. Explicitly,
\begin{equation}\label{eq:hgd.coac}
(M\ot_R \gamma_L)\circ \varrho_R = (\varrho_R \ot_L H) \circ \varrho_L
\qquad \textrm{and}\qquad
(M\ot_L \gamma_R)\circ \varrho_L = (\varrho_L \ot_R H) \circ \varrho_R.
\end{equation}
Morphisms of $\hH$-comodules are $\hH_R$-comodule maps as well as
$\hH_L$-comodule maps. The category of right $\hH$-comodules is denoted by
$\M^\hH$. 

The category ${}^\hH\M$ of left $\hH$-comodules is defined symmetrically.
\end{definition}

Note that since the right $R$- and $L$-actions on $H$ commute, also any right
$\hH$-comodule is a right $R \ot_k L$-module.  

\begin{remark}\label{rem:antip.functor}
The antipode $S$ in a Hopf algebroid $\hH$ defines a functor ${}^\hH\M\to
\M^\hH$. Indeed, if $M$ is a left $\hH$-comodule, with $\hH_R$-coaction 
$m\mapsto m^{[-1]}\ot_R m^{[0]}$ and $\hH_L$-coaction $m\mapsto m_{[-1]}\ot_L
m_{[0]}$, then it is a right $\hH$-comodule with right $R$-action
$mr:=\pi_L(t_R(r)) m $, right $L$-action $ml:=\pi_R(t_L(l))m$ and respective
coactions 
\begin{equation}\label{eq:L_R_coac}
m\mapsto m_{[0]} \ot_R S(m_{[-1]})
\qquad \textrm{and}\qquad 
m\mapsto m^{[0]}\ot_L S(m^{[-1]}).
\end{equation}
Left $\hH$-comodule maps are also right $\hH$-comodule maps for
these coactions. 

A functor $\M^\hH\to {}^\hH\M$ is constructed symmetrically. 
\end{remark}

\begin{proposition}\label{prop:hgd.coinv}
Let $\hH$ be a Hopf algebroid and $(M,\varrho_L,\varrho_R)$ be
a right $\hH$-comodule. Then any coinvariant of the $\hH_R$-comodule
$(M,\varrho_R)$ is coinvariant also for the $\hH_L$-comodule $(M,\varrho_L)$.

If moreover the antipode of $\hH$ is bijective then coinvariants of the
$\hH_R$-comodule $(M,\varrho_R)$ and the $\hH_L$-comodule $(M,\varrho_L)$
coincide. 
\end{proposition}

\begin{proof}
For a right $\hH$-comodule $(M,\varrho_L,\varrho_R)$, consider the map 
\begin{equation}\label{eq:Phi}
\Phi_M:M\ot_R H \to M\ot_ L H,\qquad m\ot_R h \mapsto \varrho_L(m)
S(h), 
\end{equation}
where $H$ is a left $L$-module via the source map $s_L$ and a left $R$-module 
via the target map $t_R$, and $M\ot_L H$ is understood to be a right
$H$-module via the second factor. 
Since $\Phi_M(\varrho_R(m))=m\ot_L 1_H$ and $\Phi_M(m\ot_R 1_H)=\varrho_L(m)$,
we have the first claim in Proposition \ref{prop:hgd.coinv} proven. 
In order to prove the second assertion, note that if $S$ is an isomorphism
then so is $\Phi_M$, with the inverse $\Phi_M^{-1}(m\ot_L h) = S^{-1}(h) 
\varrho_R(m)$, where $M\ot_R H$ is understood to be a left $H$-module via the
second factor. 
\end{proof}

Although the functors $U$ and $V$ in Theorem \ref{thm:pure} below are not 
known to exist without further assumptions, they exist in all known examples 
of Hopf algebroids and they establish
isomorphisms between the categories of $\hH_L$-comodules, $\hH_R$-comodules
and $\hH$-comodules.
In the following theorem $F_R$ and $F_L$ denote the
forgetful functors $\M^{\hH_R} \to \M_k$ and $\M^{\hH_L} \to \M_k$,
respectively, while  $G_R$ and $G_L$ denote the forgetful functors $\M^\hH\to
\M^{\hH_R}$ and $\M^\hH\to \M^{\hH_L}$, respectively.   
$H$ is regarded as an $R$-bimodule via right multiplication by $s_R$ and
$t_R$ and an $L$-bimodule via left multiplication by $s_L$ and $t_L$.

\begin{theorem}\label{thm:pure}
Consider a Hopf algebroid $\hH$.

(1) If the equaliser 
\begin{equation}\label{eq:R.M}
\xymatrix{
M \ar[rr]^-{\varrho_R}&&
M\ot_R H \ar@<2pt>[rr]^-{\varrho_R\ot_R H}\ar@<-2pt>[rr]_-{M\ot_R \gamma_R}&&
M\ot_R H \ot_R H
}
\end{equation}
in $\M_L$ is $H\ot_L H$-pure, i.e. it is preserved by the functor $(-) \ot_L
H\ot_L H:\M_L\to \M_L$, for any right $\hH_R$-comodule $(M,\varrho_R)$, then
there exists a functor $U:\M^{\hH_R}\to \M^{\hH_L}$, such that $F_L \circ U
=F_R$ and $U\circ G_R=G_L$. In particular, $G_R$ is full.

(2) If the equaliser 
$$
\xymatrix{
N \ar[rr]^-{\varrho_L}&&
N\ot_L H \ar@<2pt>[rr]^-{\varrho_L\ot_L H}\ar@<-2pt>[rr]_-{N\ot_L \gamma_L}&&
N\ot_L H \ot_L H
}
$$
in $\M_R$ is $H\ot_R H$-pure, i.e. it is preserved by the functor $(-)\ot_R
H\ot_R H:\M_R\to \M_R$, for any right $\hH_L$-comodule $(N,\varrho_L)$, then
there exists a functor $V:\M^{\hH_L}\to \M^{\hH_R}$, such that $F_R \circ V
=F_L$ and $V\circ G_L=G_R$. In particular, $G_L$ is full.

(3) If both purity assumptions in parts (1) and (2) hold, then the forgetful
functors $G_R:\M^\hH\to \M^{\hH_R}$ and $G_L:\M^\hH\to \M^{\hH_L}$
are isomorphisms, hence $U$ and $V$ are inverse isomorphisms.
\end{theorem}

\begin{proof}
(1) 
Recall that \eqref{eq:R.M} defines the $\hH_R$-cotensor product $M
\Box_{\hH_R} H\simeq M$. By (2.7), $H$ is an $\hH_R$-$\hH_L$ bicomodule, with
left coaction $\gamma_R$ and right coaction $\gamma_L$. Thus in light of
\cite[22.3]{BrzWis:cor} and its Erratum, we can define a desired functor
$U:=(-) \Box_{\hH_R} H$. Clearly, it satisfies $F_L \circ U = F_R$.
For an $\hH$-comodule $(M,\varrho_L,\varrho_R)$, the coaction on the
$\hH_L$-comodule $U\big( G_R(M,\varrho_L,\varrho_R)\big)=U(M,\varrho_R)$ is
given by  
\begin{equation}\label{eq:rho'}
\xymatrix{
M \ar[r]^-{\varrho_R}&
M\Box_{\hH_R} H \ar[rr]^-{M\Box_{\hH_R} \gamma_L}&&
M\Box_{\hH_R} (H\ot_L H) \ar[r]^-{\simeq}&
(M\Box_{\hH_R} H)\ot_L H \ar[rr]^-{M\ot_R \pi_R \ot_L H}&&
M\ot_L H,
}
\end{equation}
where in the third step we used that since the equaliser \eqref{eq:R.M} is
$H\ot_L H$-pure, it is in particular $H$-pure. Using that $\varrho_R$ is a
right $\hH_L$-comodule map and counitality of $\varrho_R$, we conclude that
\eqref{eq:rho'} is equal to $\varrho_L$. Hence $U\circ G_R =G_L$. This proves
that for any two $\hH$-comodules $M$ and $M'$, and any $\hH_R$-comodule map
$f:M\to M'$, $U(f)=f$ is an $\hH_L$-comodule map hence an $\hH$-comodule map,
i.e. that $G_R$ is full.
Part (2) is proven symmetrically.

(3) 
For the functor $U$ in part (1) and a right $\hH_R$-comodule $(M,\varrho_R)$,
denote $U(M,\varrho_R)=:(M,\varrho_L)$. 
With this notation, define a functor ${\widehat G}_R:\M^{\hH_R}\to
\M^\hH$, with object map $(M,\varrho_R)\mapsto (M,\varrho_L,\varrho_R)$, and
acting on the morphisms as the identity map. Being coassociative, $\varrho_R$
is an $\hH_R$-comodule map, so by part (1) it is an $\hH_L$-comodule
map. Symmetrically, by part (2) $\varrho_L$ is an $\hH_R$-comodule map. So
${\widehat G}_R$ is a well defined functor. We claim that it is the inverse of
$G_R$. Obviously, $G_R \circ {\widehat G}_R$ is the identity functor. In the
opposite order, 
note that by construction $G_L \circ {\widehat G}_R =U$. Therefore, $G_L \circ
{\widehat G}_R \circ G_R = U\circ G_R =G_L$, cf. part (1). That is, ${\widehat
  G}_R \circ G_R$ takes an $\hH$-comodule $(M,\varrho_L,\varrho_R)$ to the
same $\hH_L$-comodule $(M,\varrho_L)$. Since ${\widehat G}_R \circ G_R$
obviously takes $(M,\varrho_L,\varrho_R)$ to the same $\hH_R$-comodule
$(M,\varrho_R)$ as well, we conclude that also
${\widehat G}_R\circ G_R$ is the identity functor.

In a symmetrical way, in terms of the functor $V(N,\varrho_L)=:(N,\varrho_R)$
in part (2), one constructs $G_L^{-1}$ with object map $(N,\varrho_L)\mapsto
(N,\varrho_L,\varrho_R)$, and acting on the morphisms as the identity map. 
The identities $G_L \circ G_R^{-1}=U$ and $G_R \circ G_L^{-1}=V$ prove that
$U$ and $V$ are mutually inverse isomorphisms, as stated.
\end{proof}

\begin{example} We list some families of Hopf algebroids $\hH$ in which the
  purity (i.e. equaliser-preserving) conditions in Theorem \ref{thm:pure} hold.

(1) All purity conditions in Theorem \ref{thm:pure} hold if $H$ is {\em flat}
  as a left $L$- and a left $R$-module. Indeed, in this case the functors
  $(-)\ot_L H:\M_L \to \M_L$ and $(-)\ot_R H:\M_R \to \M_R$ preserve any
  equaliser. In particular, Frobenius Hopf algebroids in \cite{BalSzl:fgal}
  (being finitely generated and projective) are flat. 

(2) Weak Hopf algebras, introduced in \cite{BNSz:WHAI}, determine Hopf
  algebroids over Frobenius-separable base algebras $L\simeq R^{op}$,
  cf. \cite[4.1.2]{Bohm:HoA}. Recall that Frobenius separability of a
  $k$-algebra $R$ means the existence of a $k$-module map $\psi:R \to k$ and
  an element $\sum_i e_i \ot_k f_i \in R \ot_k R$, such that 
$$
\sum_i \psi(r e_i) f_i = r = \sum_i e_i \psi(f_i r), \quad \textrm{for all }
r\in R, \qquad \textrm{and}\qquad \sum_i e_i f_i =1_R.
$$
Note that this implies that $\sum_i r e_i \ot_k f_i = \sum_i e_i \ot_k f_i r$,
for all $r\in R$ (hence the name {\em separable}).
For a Frobenius-separable algebra $R$, any right $R$-module $X$ and left
$R$-module $Y$, the canonical epimorphism $X \ot_k Y\to X \ot_R Y$ is split by
$x\ot_R y \mapsto \sum_i x e_i \ot_k f_i y$. For the base algebras of a weak
Hopf algebra, the Frobenius-separability structure arises from the restriction
of the counit and the image of the unit element $1_H$ under the
coproduct. 

Let $H$ be a weak Hopf algebra with (weak) coproduct $\Delta: H\to H \ot_k
H$. Denote its left and right (or `target' and `source') subalgebras by $L$
and $R$, respectively. (These serve as the base algebras of the corresponding
Hopf algebroid, see \cite{Bohm:HoA}.) Note that $\Delta$ is an $R$-$L$ bimodule
map.  

Any right comodule $(M,\rho)$ of a weak Hopf algebra $H$ (as a coalgebra) 
can be equipped with a right $R$-action via $mr=m_{<0>}
\varepsilon(m_{<1>}r)$, where $\rho(m)=m_{<0>} \ot_k m_{<1>}$ and
$\varepsilon$ denotes the counit of $H$.
Moreover, any right comodule $(M,\rho)$ of $H$  yields a right
coaction of the constituent right bialgebroid, by composing $\rho$ with
the (split) epimorphism $p_M:M\ot_k H \to M\ot_R H$. 
Define a right $L$-module structure on $M$ via \eqref{eq:R.M}.

For any left $L$-module $N$, consider the following diagram (in $\M_k$). 
\begin{equation}\label{eq:wha}
\xymatrix{
M\ot_L N \ar[rr]^-{p_M\circ \rho\ot_L N}\ar[d]&&
M \ot_R H \ot_L N \ar@<2pt>[rr]^-{p_M\circ \rho\ot_R H \ot_L N}
\ar@<-2pt>[rr]_-{M\ot_R p_H\circ \Delta \ot_LN}\ar[d]&&
M\ot_R H \ot_R H \ot_L N\ar[d]\\
M\ot_k N \ar[rr]_-{\rho\ot_k N}&&
M\ot_k H \ot_k N\ar@<2pt>[rr]^-{\rho \ot_k H \ot_k N}
\ar@<-2pt>[rr]_-{M\ot_k \Delta\ot_k N}&&
M\ot_k H \ot_k H \ot_k N\ .
}
\end{equation}
The $R$-actions on $H$ are given by right multiplications by the source and
target maps and the $L$-actions on $H$ are given by left multiplications.
The vertical arrows denote the sections of the canonical epimorphisms, given
by the Frobenius-separability structures of $L$ and $R$. The diagram is
easily checked to be serially commutative (meaning commutativity with either
{\em simultaneous} choice of the upper or the lower ones of the parallel
arrows). Clearly,
$$
\xymatrix{
M \ar[r]^-{\rho}& M \ot_k H \ar@<2pt>[rr]^-{\rho \ot_k H}
\ar@<-2pt>[rr]_-{M\ot_k \Delta}&&
M\ot_k H\ot_k H
}
$$
is a split equaliser in $\M_k$ (with splitting provided by the counit
of $H$), hence the bottom row of the diagram in \eqref{eq:wha} is an
equaliser. This implies
that also the top row is an equaliser, so in particular the purity conditions
in Theorem \ref{thm:pure} (1) hold. The conditions in Theorem \ref{thm:pure}
(2) are verified by a symmetrical reasoning.

(3) For any $k$-algebra $L$, the tensor product algebra $H:=L\ot_k L^{op}$
carries a Hopf algebroid structure, see \cite[4.1.3]{Bohm:HoA}. Since the left
$L$-action on $H$ is given 
by multiplication in the first factor, the functors $F((-)\ot_L H \ot_L H)$
and $F(-)\ot_k L\ot_k L:\M_L \to \M_k$ are naturally
isomorphic, where $F:\M_L\to \M_k$ denotes the forgetful
functor. The forgetful functor $F$ has a left adjoint, hence it preserves any
equaliser. The functor $F$ takes \eqref{eq:R.M} to a split equaliser (with
splitting provided by $\pi_R$), which is then preserved by any
functor. This proves that $F(-)\ot_k L\ot_k L$
and hence $F((-)\ot_L H \ot_L H)$ preserve \eqref{eq:R.M}. Since $F$ also
reflects equalisers, we conclude that \eqref{eq:R.M} is preserved by $(-)\ot_L
H \ot_L H:\M_L \to \M_L$. The purity conditions in Theorem \ref{thm:pure}(2)
are proven to hold similarly. 

(4)
In \cite[Corrigendum]{ArdBohmMen:Sch.type}, the purity conditions in Theorem
\ref{thm:pure} are proven to hold for a Hopf algebroid whose constituent
$R$-coring (equivalently, the constituent $L$-coring) is coseparable.
\end{example}

\begin{theorem}\label{thm:hgd.com.mon}
For any Hopf algebroid ${\mathcal H}$, ${\M}^{\mathcal H}$ is a  monoidal
category. Moreover, there are strict monoidal forgetful functors rendering
commutative the following diagram:
$$
\xymatrix{
{\M}^{\hH}\ar[r]^-{G_R}\ar[d]_-{G_L}& {\M}^{\hH_R}\ar[d]\\
{\M}^{\hH_L}\ar[r]& {}_R{\M}_R \ .
}
$$
\end{theorem}

\begin{proof}
Commutativity of the diagram follows by
comparing the unique $R$-actions that make $R$-bilinear the $\hH_R$-coaction
and the $\hH_L$-coaction in an $\hH$-comodule, respectively, (see
\eqref{eq:RacM} and \eqref{eq:LacN}, and the algebra isomorphism
\eqref{eq:uv}).  
Strict monoidality of the functors on the right hand side and in the bottom
row follows by \cite[Theorem 5.6]{Schau:strthm} (and its application to 
the opposite of the bialgebroid $\hH_L$), cf. Section \ref{sec:combgd}. 
In order to see strict monoidality of the remaining two functors $G_R$ and
$G_L$, recall that by \cite[Theorem 5.6]{Schau:strthm} (applied to $\hH_R$ and
the opposite of $\hH_L$), the $R$-module tensor product of any two
$\hH$-comodules is an $\hH_R$-comodule and an $\hH_L$-comodule, via the
diagonal coactions, cf. \eqref{eq:MNcR} and \eqref{eq:L-mon}. It is
straightforward to check compatibility of these coactions in the sense of
Definition \ref{def:hgd_comod}. Similarly, $R$($\simeq L^{op}$) is known to be
an $\hH_R$-comodule and an $\hH_L$-comodule, and compatibility of the coactions
is obvious. Finally, the $R$-module tensor product of $\hH$-comodule maps is an
$\hH_R$-comodule map and an $\hH_L$-comodule map by \cite[Theorem
  5.6]{Schau:strthm}. Thus it is an $\hH$-comodule map. By Theorem
\cite[Theorem 5.6]{Schau:strthm} also the coherence natural transformations in
${}_R \M_R$ are $\hH_R$- and $\hH_L$-comodule maps, so $\hH$-comodule maps,
what completes the proof.
\end{proof}

\begin{definition}\label{def:hgd.comod.alg}
A {\em right comodule algebra} of a Hopf algebroid ${\mathcal H}$ is a
monoid in the monoidal category ${\M}^\hH$ of right
$\hH$-comodules. Explicitly, an $R$-ring $(A,\mu,\eta)$, such that $A$ is a
right $\hH$-comodule and $\eta:R \to A$ and $\mu:A\ot_R A \to A$ are right
$\hH$-comodule maps. Using the notations $a\mapsto a^{[0]}\ot_R a^{[1]}$ and
$a\mapsto a_{[0]}\ot_L a_{[1]}$ for the $\hH_R$- and $\hH_L$-coactions,
respectively, $\hH$-colinearity of $\eta$ and $\mu$ means the identities, for
all $a,a'\in A$, 
\begin{eqnarray*}
&{1_A}^{[0]} \ot_ R {1_A}^{[1]} = 1_A\ot_ R 1_H,\qquad 
&(aa')^{[0]} \ot_ R (aa')^{[1]} = a^{[0]} a^{\prime [0]} \ot_ R a^{[1]}
  a^{\prime [1]}\\
&1_{A[0]} \ot_ L 1_{A{[1]}} = 1_{A}\ot_ L 1_H,\qquad 
&(aa')_{[0]} \ot_ L (aa')_{[1]} = a_{[0]} a'_{[0]} \ot_ L a_{[1]}
  a'_{[1]} .
\end{eqnarray*}
Symmetrically, a {\em left $\hH$-comodule algebra} is a monoid in
${}^\hH\M$. 
\end{definition}

If $A$ is a right comodule algebra of a Hopf algebroid $\hH$, with
$\hH_R$-coinvariant subalgebra $B$, then we say that $B\subseteq A$ is a
(right) {\em  $\hH$-extension}. 

The functors in Remark \ref{rem:antip.functor} induced by the antipode are
checked to be 
strictly anti-monoidal. Therefore, the opposite of a right $\hH$-comodule
algebra, with coactions in Remark \ref{rem:antip.functor}, is a left
${\mathcal H}$-comodule algebra and conversely.

Whenever the antipode in a Hopf algebroid $\hH$ is bijective, it induces 
strict anti-monoidal isomorphisms ${}^\hH\M\to \M^\hH\simeq {}^{\hH_{cop}}\M$
and $\M^\hH\to {}^\hH\M \simeq \M^{\hH_{cop}}$.} 

\section{$\hH$-cleft extensions}\label{sec:Hcleft}

Recall that to an $L$-ring $A$ (with multiplication $\mu: A\ot_L
A\to A$ and unit map $\eta:L\to A$) and an $L$-coring $H$ (with
comultiplication $\gamma:H\to H\ot_L H$ and counit $\pi:H\to L$),
one associates a {\em convolution algebra} ${\rm Hom}_{L,L}(H,A)$,
with multiplication $j\diamond j'\colon =\mu\circ (j\sstac{L}
j')\circ \gamma$ and unit $\eta\circ\pi$. The first aim of this
section is to develop a generalisation of the notion of a
convolution algebra and, in particular, of a convolution inverse
suitable for Hopf algebroids. 
This will make it possible to interpret in particular the antipode of a Hopf
algebroid as the convolution inverse  of the identity map.

As explained in Section~\ref{sec:hgd}, a Hopf algebroid is built on
a $k$-module with two coring structures. Although we are
primarily interested in  Hopf algebroids, in general there is no need
to put any special restrictions on these coring structures. Dually, one can
consider
a $k$-module with ring structures over two different rings. In this more
general situation the convolution algebra (which is simply a $k$-linear
category with a single object) can be generalised to a Morita context (i.e. a
$k$-linear category with two objects).
The notion of a {\em convolution inverse} is introduced within this
{\em convolution category}.

Let $L$ and $R$ be $k$-algebras and  let $H$ and $A$ be $k$-modules.
Assume that $A$ is an $L$-ring  (with multiplication $\mu_L: A\ot_L A\to A$
and unit  $\eta_L:L\to A$) and an $R$-ring  (with multiplication
$\mu_R:A\ot_R A\to A$ and unit $\eta_R:R\to A$).
Assume that $A$ is an $L$-$R$ and $R$-$L$ bimodule
with respect to the corresponding module structures,  $\mu_L$
 is  $R$-$R$ bilinear and $\mu_R$ is  $L$-$L$ bilinear,
and that
\begin{equation}\label{eq.prod.ass2}
\mu_L\circ (A\stac{L} \mu_R)= \mu_R\circ (\mu_L\stac{R} A), \qquad
\mu_R\circ (A\stac{R} \mu_L)= \mu_L\circ (\mu_R\stac{L} A).
\end{equation}
Dually, assume that $H$ is an $L$-coring  (with comultiplication
$\gamma_L:H\to H\ot_L H$ and counit $\pi_L:H\to L$) and an $R$-coring
(with
comultiplication
$\gamma_R:H\to H\ot_R H$ and counit $\pi_R:H\to R$). Assume
further that $H$ is an $L$-$R$ and $R$-$L$ bimodule with respect to the
corresponding module structures, such that $\gamma_L$ is $R$-$R$ bilinear,
$\gamma_R$ is $L$-$L$ bilinear and
\begin{equation}\label{eq.coprod.ass2}
(H\stac{L} \gamma_R)\circ \gamma_L = (\gamma_L\stac{R} H)\circ \gamma_R,
\qquad
(H\stac{R} \gamma_L)\circ \gamma_R = (\gamma_R\stac{L} H)\circ \gamma_L.
\end{equation}
To the above data one associates a
$k$-linear {\em convolution category} $\Conv H A$ as
follows.
 $\Conv H A$ has two objects, $R$ and $L$,
 and morphisms
 $$
 \Conv H A(P,Q) = \lrhom QPHA, \qquad P,Q\in \{L,R\},
 $$
  with composition $\diamond$, defined for all
  $\phi \in \lrhom PQHA $ and  $\psi \in \lrhom QSHA$,  $ P,Q,S \in \{L,R\}$,
  $$
  \phi\diamond\psi = \mu_Q\circ (\phi\stac{Q} \psi)\circ \gamma_Q
    \in \lrhom PSHA.
  $$
Note that the identity morphism in $\Conv H A(P,P)$ is $\eta_P\circ\pi_P$. The
conditions \eqref{eq.prod.ass2} and \eqref{eq.coprod.ass2} together with
coassociativity of the coproducts in $H$ and associativity of products in $A$
ensure that
the composition $\diamond$ is an associative operation.
\begin{definition}
\label{def.conv.inv}
Let $\Conv H A$ be a convolution category
and let $j$ be a morphism
in $\Conv H A$.  A retraction of $j$ in $\Conv H A$ is called
 a {\em left convolution inverse of $j$} and a section of $j$ in $\Conv H A$
 is called
 a {\em right convolution inverse of $j$}.
 If $j$ is an isomorphism in $\Conv H A$,
 then it is said to be {\em convolution invertible}; its inverse is called the
 {\em convolution inverse of $j$} and is denoted  by $\conv{j}$.
 \end{definition}
\begin{remark}(1) A $k$-linear category with a single object $a$ can be
  identified 
  with the $k$-algebra ${\rm End}(a)$ of the morphisms in the category. In a
  similar manner, a $k$-linear category with two objects $a$ and $b$ can be
  identified with a Morita context as follows. The composition of
  morphisms makes  $k$-modules ${\rm Hom}(a,b)$ and ${\rm Hom}(b,a)$
  bimodules for  $k$-algebras ${\rm End}(a)$ and ${\rm
  End}(b)$. Furthermore, the restriction of the composition to the map
  $$
  {\rm Hom}(a,b)\stac{k} {\rm Hom}(b,a)\to {\rm End}(b)
  $$
  is an ${\rm End}(a)$-balanced ${\rm End}(b)$-bimodule map.
  That is, it is a composite of the canonical epimorphism
 $
 {\rm Hom}(a,b)\sstac{k} {\rm Hom}(b,a)\to {\rm Hom}(a,b)\sstac{{\rm End}(a)}
  {\rm Hom}(b,a),
  $
  and an ${\rm End}(b)$ -bimodule map, $F_a:{\rm
  Hom}(a,b)\sstac{{\rm End}(a)} {\rm Hom}(b,a)\to {\rm End}(b)$. Similarly,
  the map
  $$
  {\rm Hom}(b,a)\stac{k}{\rm Hom}(a,b)\to {\rm End}(a),
  $$
  obtained
  by restricting the composition, factors through the
  canonical epimorphism and the ${\rm End}(a)$-bimodule map, $F_b:{\rm
  Hom}(b,a)\sstac{{\rm End}(b)}{\rm   Hom}(a,b)\to {\rm End}(a)$.
  Using the associativity of the composition of
  the morphisms in a category, one easily checks that the 6-tuple $({\rm
  End}(a), {\rm End}(b),{\rm Hom}(b,a),$ ${\rm Hom}(a,b), F_a,F_b)$ is a
  Morita context. Clearly, there is a category of this kind behind any Morita
  context.

  In particular, the convolution category $\Conv H A$ can be identified with a
  Morita context connecting convolution algebras ${\rm Hom}_{L,L}(H,A)$
  and ${\rm Hom}_{R,R}(H,A)$.

(2) In the case $L=R$, $\gamma_L=\gamma_R$, $\pi_L=\pi_R$,
$\mu_R=\mu_L$, $\eta_R=\eta_L$, i.e.\ when there is one, say, $L$-coring
$H$ and one, say, $L$-ring $A$, the convolution category
$\Conv HA$ consists of a single object. The algebras in the corresponding
Morita context are both equal to the convolution algebra
${\rm Hom}_{L,L}(H,A)$, the bimodules are the regular bimodules and the
connecting homomorphisms are both equal to the identity map of ${\rm
Hom}_{L,L}(H,A)$. In a word: the Morita context reduces to the convolution
algebra. Thus an $L$-$L$ bimodule map $j$
is convolution invertible in the sense of Definition~\ref{def.conv.inv} if and
only
if it is an invertible element of the convolution algebra $\lrhom LLHA$.

(3) Conditions \eqref{eq.prod.ass2}, imposed on the $R$-ring and
$L$-ring structures of $A$, imply that the underlying $k$-algebras
are isomorphic via the map $A\ni a\mapsto \mu_R(a\sstac{R}
\eta_L(1_L))$, with the inverse $a\mapsto \mu_L(a\sstac{L}
\eta_R(1_R))$.

(4) Conditions \eqref{eq.coprod.ass2}, imposed on the two coring structures
of $H$, imply that the
$L$-coring $H$ is a left (and right) extension of the $R$-coring $H$, while the
$R$-coring $H$ is a right (and left) extension of the $L$-coring
$H$
in the sense of \cite{Brz:not}, with the coactions given by the coproducts.
\end{remark}

We can now exemplify the contents of Definition~\ref{def.conv.inv}
with the main case of interest, 
whereby the coring structures on $H$ constitute bialgebroids.
Consider a right
bialgebroid $\hH_R=(H,R,s_R,t_R,$ $\gamma_R,\pi_R)$ and a left bialgebroid
$\hH_L=(H,L,s_L,t_L,\gamma_L,\pi_L)$ on the same total algebra $H$, which
satisfy conditions \eqref{ax:st} and \eqref{ax:cpcolin}. 
In this situation, compatibility conditions 
for coring structures on $H$ in the definition of a convolution category,
including equations \eqref{eq.coprod.ass2}, are satisfied. For a target
of convolution invertible maps take  an $R\ot_k
L$-ring $A$. In this case the unit maps $\eta_R$ and ${\eta}_L$ are
obtained
as the restrictions of the unit map $R\ot_k L\to A$ to $R\ot_k
1_L$ and to $1_R\ot_k L$, respectively. There is no need to
 distinguish between the products of $A$ as an $R$-ring and as an
 $L$-ring, so we write simply $\mu_A$ for the product in $A$, and it becomes
 clear from the context, how this should be understood.
 Since we are dealing with a single product, it makes sense
 to denote the action of $\mu_A$ on elements by juxtaposition.
 One immediately checks
 that all the compatibility conditions between the $L$-, and
$R$-ring structures
 on $A$
 in the definition of a convolution category
  are satisfied, in particular \eqref{eq.prod.ass2}
 follow by the associativity of $\mu_A$. 
All this means that, for 
two bialgebroids $\hH_L$ and $\hH_R$ on the same total algebra $H$, such that
\eqref{ax:st} and \eqref{ax:cpcolin} hold, 
and an $R\ot_k L$-ring $A$,  there is a convolution category $\Conv HA$.
 We can now make explicit the contents of Definition~\ref{def.conv.inv} in this
 case. This essentially means describing explicitly all the 
 $L$ and $R$ actions involved.

For left and right bialgebroids $\hH_L$ and $\hH_R$ on the same total
algebra $H$, such 
that \eqref{ax:st} and \eqref{ax:cpcolin} hold, 
and an $R\ot_k
L$-ring $A$, a  map $j:H\to A$ is an $L$-$R$ bimodule map provided
\begin{equation}
j\left(s_L(l) \,h\,s_R(r)\right)={\eta}_L(l)\,j(h)\,\eta_R(r),\qquad
\textrm{for\  all\ } l\in L, r\in R,\  h\in H.\label{eq:js}
\end{equation}
Similarly, ${\tilde j}:H\to A$ is an $R$-$L$ bimodule map if
\begin{equation}
 {\tilde j}\left(t_L(l) \,h\, t_R(r)\right)={\eta}_R(r)
{\tilde j}(h){\eta}_L(l),\qquad  \ \
\textrm{for\  all\ } l\in L, r\in R,\  h\in H.\label{eq:tjt}
\end{equation}
A right convolution inverse of  $j\in \lrhom LRHA$ is a map
${\tilde j}\in \lrhom RLHA$ such that
\begin{equation}\label{eq:rinv}
\mu_A\circ (j\stac{R} {\tilde j})\circ\gamma_R=\eta_L\circ \pi_L.
\end{equation}
A  left convolution inverse of $j$ is a map
${\hat j}\in \lrhom RLHA$ such that
\begin{equation}\label{eq:linv}
\mu_A\circ ({\hat j}\stac{L} j)\circ\gamma_L=\eta_R\circ \pi_R.
\end{equation}
Obviously, by the associativity of the composition in $\Conv HA$,
if a map $j:H\to A$ satisfying (\ref{eq:js}) has both left and
right convolution inverses, then they coincide and hence the
convolution inverse of an $L$-$R$ bimodule map $j$ is unique.

\begin{example}\label{ex:hgdax}
Let $\hH_L=(H,L,s_L,t_L,\gamma_L,\pi_L)$ be a left bialgebroid and
$\hH_R=(H,R,s_R,t_R,\gamma_R,$ $\pi_R)$ be a right bialgebroid, on the same
total algebra $H$. Assume that the compatibility conditions \eqref{ax:st} and 
\eqref{ax:cpcolin} hold. Consider the $R\sstac{k} L$-ring structure on $H$,
defined by the unit map $R\sstac{k} L\to H$, $r\sstac{k} l\mapsto
  s_R(r)s_L(l)\equiv s_L(l) s_R(r)$. It gives rise to a convolution category 
$\Conv HH$. In light of \eqref{eq:js}, the identity map of $H$ is an element
of $\Conv HH(R,L)$.  By \eqref{eq:tjt}, \eqref{eq:rinv} and \eqref{eq:linv},
the identity map possesses a  
convolution inverse $S$ if and only if the first equality in \eqref{ax:Sbilin} 
and \eqref{ax:antip} hold true. Hence it follows by Remark \ref{rem:axred}
that $(\hH_L,\hH_R,S)$ is a Hopf algebroid if and only if $S$ is convolution
inverse of the identity map in $\Conv HH$ (in the same way as the antipode 
of a Hopf algebra $H$ over a ring $k$ is the inverse of the identity map in
the convolution algebra $\mathrm{End}_k(H)$).
\end{example}

\begin{example} \label{ex:smash}
Example \ref{ex:hgdax} can be extended as follows.
Take a  Hopf algebroid $\hH = (\hH_L,\hH_R,S)$ and a left
$\hH_L$-module algebra $B$. 
(The role of $B$ is played by the base algebra $L$ in Example
  \ref{ex:hgdax}.)
The smash product algebra $A\colon =B
\rtimes H$ is defined as the $k$-module $B\sstac{L} H$ with
product
$$
(b\rtimes h)(b'\rtimes h')\colon = b\,(h\di\cdot b')\rtimes h\dii h',
$$
(cf. \cite[Section 2.3]{KadSzl:D2bgd}). Here, the left $L$-module
structure on $H$ is given by the multiplication by $s_L(l)$ on the
left. $A$ is an $R$-ring with ${\eta}_R(r)=1_B\rtimes s_R(r)$ (and
hence an $L^{op}$-ring with unit $l\mapsto 1_B\rtimes t_L(l)$) and
an $L$-ring with unit $\eta_L(l)=1_B\rtimes s_L(l)$. Since the
elements $\eta_R(r)$ and $\eta_L(l)$ commute in $A$, for any $r\in
R$ and $l\in L$, $A$ is an $R\sstac{k}L$-ring.

The $L$-$R$ bimodule map $j:H\to A$, $h\mapsto 1_B\rtimes
  h$ is convolution invertible with the inverse $\conv{j}:H\to A$,
$h\mapsto 1_B\rtimes S(h)$ (cf. \eqref{ax:antip}).
\end{example}

The notion of a convolution inverse, once established, plays the
fundamental role in the definition of a cleft extension of a Hopf algebroid,
which we describe presently.
Let $\hH$ be a Hopf algebroid and $A$ a right
$\hH$-comodule algebra. Then $A$ is, in particular, an $R$-ring. The unit of
this $R$-ring, the algebra homomorphism $R\to A$, is
denoted by $\eta_R$.
{\color{blue}
The $k$-subalgebra of $\hH_R$-coinvariants 
(whose elements are also $\hH_L$-coinvariants, cf. Proposition
\ref{prop:hgd.coinv}) in $A$ is denoted by $B$.} 

Assume that $A$ is also an $L$-ring, with unit
${\eta}_L:L \to A$, and $B$ is an $L$-subring of $A$. The latter implies that
both
the $\hH_R$-coaction $\roA$, and the $\hH_L$-coaction $\lambda^A$ are left
$L$-linear. Since $\roA$ is $R$-$R$ bilinear (cf. \eqref{eq:robilin}),
$$
\roA(b\eta_R(r))=b\stac{R} s_R(r)=\roA(\eta_R(r) b),\qquad
\textrm{for\ all\ }r\in R,\ b\in B.
$$
Thus it follows that $B$ is in the commutant of the image of $\eta_R$.

Recall from Section \ref{sec:combgd} that any right $\hH$-colinear map
$j:H\to A$
is right $R$-linear in the sense of (\ref{eq:js}) and left $R$-linear in the
sense that
\begin{equation}\label{eq:fourth}
j(s_R(r) h)=\eta_R(r) j(h),
\end{equation}
for all $r\in R$ and $h\in H$ (cf. \eqref{eq:RacM}).
\begin{definition} \label{def:cleft}
Let $\hH$ be a Hopf algebroid and $A$ a right
$\hH$-comodule algebra. Denote by $\eta_R(r)=r\cdot 1_A=1_A\cdot r$ the unit
map of the corresponding $R$-ring structure of $A$. Let $B$ be the
subalgebra of 
{\color{blue}
$\hH_R$}-coinvariants in $A$. The extension $B\subseteq A$ is
called {\em $\hH$-cleft} if
\begin{blist}
\item $A$ is an $L$-ring
   (with unit ${\eta}_L:L \to A$) and $B$ is an $L$-subring of $A$;
\item there exists a convolution invertible left $L$-linear right
  $\hH$-colinear morphism $j:H\to A$.
\end{blist}
A map $j$, satisfying condition (b), is called a {\em cleaving map}.
\end{definition}
Condition (b) in Definition \ref{def:cleft} means, in particular,
that a cleaving map is $L$-$R$ bilinear in the sense of
\eqref{eq:js}.
\begin{example}
Consider a smash product algebra $A=B\rtimes H$ of Example
\ref{ex:smash}. Similarly to \cite[Example 3.7]{Bohm:gal}, $A$ is
a right $\hH$-comodule algebra with $\hH_R$-coaction $B\sstac{L}
\gamma_R$
{\color{blue}
and $\hH_L$-coaction $B\sstac{L} \gamma_L$}. 
The 
{\color{blue}
coinciding subalgebra of $\hH_R$-coinvariants and $\hH_L$-coinvariants} 
in $A$ is $B\rtimes 1_H$. It is an $L$-subring of $A$. Since the convolution
invertible map $j:H\to A$, $h\mapsto 1_B\rtimes h$ in Example
\ref{ex:smash} is right $\hH$-colinear, $B\subseteq A$ is an
$\hH$-cleft extension.

In particular, let $N\subseteq M$ be a depth 2 (or D2, for short) extension of
algebras \cite[Definition~3.1]{KadSzl:D2bgd}.
It has been proven in \cite[Theorem 4.1]{KadSzl:D2bgd} that
the algebra $\lrend NNM$ of $N$-$N$ bilinear endomorphisms of $M$
is a left bialgebroid  and $M$ is its left module algebra.
By \cite[Corollary 4.5]{KadSzl:D2bgd},
the algebra $\rend NM$ of right $N$-linear endomorphisms of $M$,
with multiplication given by composition, is isomorphic to the smash
product algebra $M \rtimes \lrend NNM$.

If the D2 extension
$N\subseteq M$ is also a Frobenius extension then  $\lrend NNM$ is a
Hopf algebroid. Hence we can conclude that for any D2 Frobenius
extension $N\subseteq M$, the extension $M \subseteq \rend NM$
(where the inclusion is given by the left multiplication) is a cleft extension
of the Hopf algebroid $\lrend NNM$.
\end{example}

\begin{lemma}\label{lem:third}
Let $\hH$ be a Hopf algebroid and $B\subseteq A$ an $\hH$-cleft
extension, with a cleaving map $j$.
Then
\begin{equation}
\conv{j}(t_R(r)\, h)= \conv{j}(h)\, \eta_R(r),\qquad
\textrm{for\ all\ } r\in R,\ h\in H. \label{eq:third}
\end{equation}
\end{lemma}
\begin{proof}
Use the counit property of $\pi_L$ (in the first equality),
right $L$-linearity of $\conv{j}$, i.e.\ \eqref{eq:tjt} (in the second one),
the fact that, since $B$ is an $L$-subring of $A$, the images of
 $\eta_L$ and of $\eta_R$ commute in $A$ (in the third one),
the assumption that $\conv{j}$ is right convolution inverse of $j$,
 i.e.\ \eqref{eq:rinv} (in the fourth one),
the left $R$ linearity of $j$, i.e. \eqref{eq:fourth} (in the
fifth one), axiom \eqref{ax:cpcolin} (in the sixth one), the
identity $\gamma_L(s_R(r)h)=h\di\sstac{L} s_R(r)h\dii$, for $h\in
 H$ and $r\in R$, and the assumption that $\conv{j}$ is a left
 convolution inverse of $j$, i.e.\ \eqref{eq:linv} (in the seventh one),
the left $R$-linearity of $\conv{j}$, i.e.\ \eqref{eq:tjt} (in the
 penultimate one)
and the counit property of $\pi_R$ (in the last one)
to compute
\begin{eqnarray*}
\conv{j}(h)\eta_R(r)
&=& \conv{j}\big( t_L(\pi_L(h\dii))h\di\big)\eta_R(r)
=\conv{j}(h\di)\eta_L(\pi_L(h\dii))\eta_R(r)\\
&=&\conv{j}(h\di)\eta_R(r)\eta_L(\pi_L(h\dii))
=\conv{j}(h\di)\eta_R(r)j({h\dii}\ui)\conv{j}({h\dii}\uii)\\
&=& \conv{j}(h\di)j\big( s_R(r) {h\dii}\ui\big)\conv{j}({h\dii}\uii)
=\conv{j}({h\ui}\di)j\big( s_R(r) {h\ui}\dii \big)\conv{j}({h}\uii)\\
&=& \eta_R\big(\pi_R( s_R(r)h\ui)\big)\conv{j}({h}\uii)
=\conv{j}\big(h\uii t_R(\pi_R(s_R(r)h\ui))\big)
=\conv{j}(t_R(r) h),
\end{eqnarray*}
for $h\in H$ and $r\in R$.
\end{proof}

In the case of  a Hopf algebra cleft extension, the convolution
inverse of a cleaving map is a right colinear map, where the right coaction
in the Hopf algebra is given by the coproduct followed by the antipode and a
flip.
In the case of a Hopf algebroid there are two coactions, one
for each constituent bialgebroid, related by 
{\color{blue}
\eqref{eq:hgd.coac}}.
The following lemma shows the behaviour
of the convolution inverse of a cleaving map  with respect
to these right coactions.
\begin{lemma} \label{lem:roaj}
Let $\hH$ be a Hopf algebroid and $B\subseteq A$ an $\hH$-cleft extension with
a cleaving map $j$.
Then, for all $h\in H$,
\begin{equation}
\roA\big( \conv{j}(h)\big)=\conv{j}(h_{(2)})\stac{R} S(h_{(1)}),
\label{eq:roaj}
\end{equation}
and
\begin{equation}
\lambda^A \big(\conv{j}(h))=\conv{j}(h\uii)\stac{L}
S(h\ui).\label{eq:lamaj}
\end{equation}
\end{lemma}
\begin{proof}
Combining the module map property of the antipode, $S(t_L(l) h)=S(h) s_L(l)$,
for all $l\in L$, $h\in H$,
with the Hopf algebroid axiom $s_L = t_R\circ\pi_R\circ s_L$ and
using (\ref{eq:third}), one shows that
the expression on the right hand side of (\ref{eq:roaj})
belongs to the appropriate $R$-module tensor product. Next using
(\ref{eq:tjt})
one finds that it is an element of the Takeuchi product $A\times_R H$, defined
in \eqref{eq:MTak}, i.e.
$$
\eta_R(r) \conv{j}(h_{(2)})\stac{R} S(h_{(1)})=
\conv{j}(h_{(2)})\stac{R} t_R(r) S(h_{(1)}),
$$
for all $r\in R$, $h\in H$.
 $A\times_R H$ is an $R\ot L$-ring with factorwise multiplication and unit maps
\begin{eqnarray*}
\eta^{\times}_R:R\to A\times_R H,\quad r\mapsto 1_A\stac{R} s_R(r)
\quad{\textrm{and}}\quad
\eta^{\times}_L:L\to A\times_R H\quad l\mapsto \eta_L(l)\stac{R} 1_H,
\end{eqnarray*}
such that $\roA:A\to A\times_R H$ is a homomorphism of $R\ot_k
L$-rings. Furthermore, $\roA\circ \conv{j}:H\to A\times_R H$ is
the convolution inverse of $\roA\circ j$. We claim that the map
$$
{\tilde \Theta}: H\to A\times_R H,\qquad h\mapsto \conv{j}(h_{(2)})\stac{R}
S(h_{(1)}) ,
$$
is a right convolution inverse of $\roA\circ j$.

Take any  $h\in H$, $l\in L$ and $r\in R$. By the Hopf algebroid
identity $\gamma_L(t_L(l) h t_R(r))=h_{(1)}t_R(r)\sstac{L} t_L(l)
h_{(2)}$, the module map property  of the antipode,
$S(ht_R(r))=s_R(r) S(h)$, and the right $L$-linearity of
$\conv{j}$ it follows that
\begin{eqnarray*}
{\tilde \Theta}(t_L(l) h t_R(r))&=&
\conv{j}(t_L(l) h_{(2)})\stac{R} S(h_{(1)} t_R(r))\\
&=& \conv{j}(h_{(2)})\eta_L(l) \stac{R} s_R(r)S(h_{(1)})
=\eta^{\times}_R(r){\tilde \Theta}(h)\eta_L(l)^{\times},
\end{eqnarray*}
that is, ${\tilde \Theta}$ satisfies (\ref{eq:tjt}).
Using the right $\hH_R$-colinearity of $j$
and the coassociativity of $\gamma_R$, one computes,
\begin{eqnarray*}
[\mu_{A\times_R H}\circ (\roA\!\circ\! j\,\,\stac{R}\,{\,\tilde \Theta})\circ
  \gamma_R](h)
&=& j({h\ui})\conv{j}({{h\uii}\uii}\dii)\stac{R} {h\uii}\ui
S({{h\uii}\uii}\di)\\
&=& j({h\ui})\conv{j}({h\uii}\dii)\stac{R} {{h\uii}\di}\ui
S({{h\uii}\di}\uii)\\
&=& j({h\ui})\conv{j}({h\uii}\dii)\stac{R} s_L\big(\pi_L({h\uii}\di)\big)\\
&=& j({h\ui})\conv{j}({h\uii}\dii)\eta_R( \pi_R(
  s_L(\pi_L({h\uii}\di))))\stac{R} 1_H\\
&=& j({h\ui})\conv{j}(  s_L(\pi_L({h\uii}\di)) {h\uii}\dii)\stac{R} 1_H\\
&=& j({h\ui})\conv{j}({h\uii})\stac{R} 1_H=\eta_L( \pi_L(h))\stac{R} 1_H=
[\eta^{\times}_L\circ \pi_L](h),
\end{eqnarray*}
where the second equality follows by
the Hopf algebroid axiom \eqref{ax:cpcolin},
the third one by the antipode axiom \eqref{ax:antip}, the fourth one
by the axiom $s_L= t_R\circ \pi_R\circ s_L$ in \eqref{ax:st},
the fifth one by (\ref{eq:third})
and the penultimate one by (\ref{eq:rinv}).
Note that (although the counit $\pi_L$ is not left $R$-linear), the
restriction of $\pi_L\ot_R H$ to the Takeuchi product $H\times_L H$ is left
$R$-linear. Hence all expressions in the above computation are meaningful.
This proves that ${\tilde \Theta}$ satisfies (\ref{eq:rinv}), hence
${\tilde \Theta}$ is a right convolution inverse of $\roA\circ j$. In view of
the uniqueness of a convolution inverse
this implies
(\ref{eq:roaj}).

{\color{blue}
$\hH_L$-colinearity (i.e. property \eqref{eq:lamaj}) of $\conv{j}$ is checked
  by similar steps:
The map $R \ot_k L\to A \times_L H$, $r\ot_k l\mapsto \eta_L(l) \ot_L s_R(r)$
equips $A \times_L H$ with an $R \ot_k L$-ring structure, such that
$\lambda^A:A \to A \times_L H$ is a morphism of $R \ot_k L$-rings. Moreover,
$\lambda^A \circ 
\conv{j}$ is the convolution inverse of $\lambda^A\circ j : H \to A \times_L
H$, and the map $h \mapsto \conv{j}(h^{(2)}) \sstac L S(h^{(1)})$ is checked
to be the left convolution inverse of $\lambda^A\circ j$. Thus the uniqueness
of a convolution inverse implies \eqref{eq:lamaj}.}
\end{proof}

\begin{remark}\label{rem:lefthand}
Recall from Section \ref{sec:comhgd} that if the antipode of a Hopf algebroid
$\hH$ is bijective, then
there exists an anti-monoidal isomorphism between the categories of right
$\hH$-comodules and right $\hH_{cop}$-comodules. Hence in this case, in light
of the explicit form \eqref{eq:L_R_coac} of the relation between the $\hH_L$
and $(\hH_R)_{cop}$-coactions, an algebra extension $B\subseteq A$ is a right
$\hH$-extension if and only if $B^{op}\subseteq A^{op}$ is a right
$\hH_{cop}$-extension. Furthermore, $B\subseteq A$ is an
$\hH$-cleft extension
if and only if $B^{op}\subseteq A^{op}$ is an $\hH_{cop}$-cleft extension.
Indeed, by Lemma \ref{lem:roaj}, if $j:H\to A$ is a cleaving map for the
$\hH$-cleft extension $B\subseteq A$, then its convolution inverse $\conv{j}$
is a cleaving map for the $\hH_{cop}$-cleft extension $B^{op}\subseteq
A^{op}$.
\end{remark}
Our next aim is to prove that, in parallel to the Hopf algebra case,
 an $\hH$-cleft extension can be equivalently characterised as a
Galois extension with the normal basis property. This is
the main result of this section. The main difference with the
Hopf algebra case is that a cleft $\hH$-extension is a Galois extension with
 respect to the {\em right} bialgebroid $\hH_R$
 but it has a normal basis property with respect
 to the base algebra $L$ of the
{\em left} bialgebroid $\hH_L$.
In preparation for this we state the following two lemmas.
\begin{lemma} \label{lem:coinv}
Let $\hH$ be a Hopf algebroid and $B\subseteq A$ an $\hH$-cleft extension
with a cleaving map $j$.
Then, for all $a\in A$,  $a^{[0]} \conv{j}(a^{[1]})\in B$.
\end{lemma}
\begin{proof}
This is checked by applying $\roA$ to $a^{[0]} \conv{j}(a^{[1]})$,
noting that $\roA$ is an algebra map and $\conv{j}$ satisfies
 equation~\eqref{eq:roaj}, and then
repeating the same steps as in the verification that
$\tilde{\Theta}$ satisfies equation~\eqref{eq:rinv}
 in the proof of Lemma~\ref{lem:roaj}.
\end{proof}
\begin{lemma} \label{lem:dirsum}
Let $\hH$ be a Hopf algebroid and $B\subseteq A$ an $\hH$-cleft extension.
Then the inclusion $B\subseteq A$ splits in the category of left $B$-modules.
If, in addition, the antipode of $\hH$ is bijective, then the inclusion
$B\subseteq A$ splits also in the category of right $B$-modules.
\end{lemma}
\begin{proof} A left $B$-linear splitting of the inclusion $B\to A$ is given
by the map
\begin{equation}\label{eq:split}
A\to B,\qquad a\ \mapsto\ a\tz\conv{j}(a\ti)j(1_H),
\end{equation}
where $j$ is a cleaving map.
The element $a\tz\conv{j}(a\ti)$ belongs to $B$ for any $a\in A$ by Lemma
\ref{lem:coinv}
and $j(1_H)$ is an element of $B$ by the colinearity of $j$ and the unitality
of $\roA$.
The left $B$-linearity of the map (\ref{eq:split}) follows by the left
$B$-linearity of $\roA$. Finally, for all $b\in B$,
$$
b\tz\conv{j}(b\ti)j(1_H)
=b\ \conv{j}(1_H)j(1_H)
=b\ \eta_R\big( \pi_R(1_H)\big)=b,
$$
where the penultimate equality follows by the fact that $\conv{j}$ is the
convolution inverse of $j$ and the unitality of $\gamma_L$.

If the antipode is bijective, then, by Remark \ref{rem:lefthand}, the map
$$
A\to B,\qquad a\ \mapsto\ \conv{j}(1_H) j\big(S^{-1}(a\bi)\big)a\bz ,
$$
is a right $B$-linear section of the inclusion $B\subseteq A$.
\end{proof}

Notice that $B\sstac{L} H$ is a left $B$-module via the regular
$B$-module structure of the first tensor factor, and -- since the
coproducts $\gamma_R$ and $\gamma_L$ are left $L$-linear -- also a
right $\hH$-comodule via the regular $\hH$- comodule structure of
the second tensor factor.

\begin{theorem}\label{thm.main.cleft}
Let $\hH$ be a Hopf algebroid and $B\subseteq A$ a right $\hH$-extension.
Then the following statements are equivalent:
\begin{zlist}
\item $B\subseteq A$ is an $\hH$-cleft extension.
\item
\begin{blist}
\item The extension $B\subseteq A$ is $\hH_R$-Galois; \item
$A\simeq B\sstac{L} H$ as left $B$-modules and
  right $\hH$-comodules.
\end{blist}
\end{zlist}
\end{theorem}
\begin{proof}
(1) $\Ra$ (2)(a) Suppose that $B\subseteq A$ is a cleft $\hH$-extension
with a cleaving map $j$ and consider the map
$$
\Phi:A\stac{R} H\to A\stac{B} A,\qquad
a\stac{R} h\mapsto a\conv{j}(h\di)\stac{B} j(h\dii).
$$
By (\ref{eq:tjt}) and \eqref{eq:js}, $\conv{j}(h\di)\sstac{L}
j(h\dii)$ is a well defined element of $A\sstac{L} A$. Since $B$
is an $L$-ring, $\Phi$ is a well defined map. We claim that $\Phi$
is the inverse of the $\hH_R$-canonical map \eqref{eq:Rcan}. Take
any $a\sstac{R} h\in A\sstac{R} H$ and compute
\begin{eqnarray*}
{\rm can}_R\big(\Phi(a\stac{R} h)\big)
&=& a\conv{j}(h\di)j(h\dii)^{[0]}\stac{R} j(h\dii)^{[1]}
= a\conv{j}(h\di)j({h\dii}\ui)\stac{R} {h\dii}\uii\\
&=& a\conv{j}({h\ui}\di)j({h\ui}\dii)\stac{R} {h}\uii
= a\eta_R\big(\pi_R(h\ui)\big)\stac{R} h\uii\\
&=& a\stac{R} h\uii t_R\big(\pi_R(h\ui)\big)=a\stac{R} h,
\end{eqnarray*}
where the second equality follows by the right $\hH_R$-colinearity
of $j$, the third one by \eqref{ax:cpcolin}, and the fourth one by
(\ref{eq:linv}). On the other hand, for all $a\sstac{B} a'\in
A\sstac{B} A$,
\begin{eqnarray*}
\Phi\big( {\rm can}_R(a\stac{B} a')\big)
&=& a a^{\prime [0]} \conv{j}({a^{\prime[1]}}\di)\stac{B}
  j({a^{\prime[1]}}\dii)
= a {{a^{\prime}}_{[0]}}^{[0]} \conv{j}({{a^{\prime}}_{[0]}}^{[1]})\stac{B}
  j({{a^{\prime}}_{[1]}})\\
&=& a \stac{B} {{a^{\prime}}_{[0]}}^{[0]} \conv{j}({{a^{\prime}}_{[0]}}^{[1]})
  j({{a^{\prime}}_{[1]}})
= a \stac{B} {{a^{\prime}}}}^{[0]} \conv{j}({{a^{\prime}}^{[1]}}_{(1)})
  j({{{a^{\prime}}^{[1]}}\dii)\\
&=& a \stac{B} {{a^{\prime}}}^{[0]} \eta_R
\big(\pi_R({{a^{\prime}}^{[1]}})\big)=
a\stac{B}a',
\end{eqnarray*}
where the second and the fourth equalities follow by the right
$\hH_L$-colinearity of $\roA$, the third one by Lemma \ref{lem:coinv}, the
fifth one by (\ref{eq:linv}) and the last one by the counitality of
$\roA$. Thus $\Phi$
is the inverse of the canonical map, as claimed.

(1) $\Ra$ (2)(b) Given a cleaving map $j$, consider the left $B$-linear
map
\begin{equation}\label{eq:kappa}
\kappa: A\to B\stac{L} H,\qquad a\mapsto
{a_{[0]}}^{[0]}\conv{j}({a_{[0]}}^{[1]})\stac{L} a_{[1]}=
a^{[0]}\conv{j}({a^{[1]}}_{(1)})\stac{L} {a^{[1]}}_{(2)}.
\end{equation}
The equality of two forms of $\kappa$ follows by the right
$\hH_L$-colinearity of $\roA$. Furthermore, Lemma \ref{lem:coinv}
implies that the image of $\kappa$ is in $B\sstac{L} H$. In the
opposite direction, define the left $B$-linear map $
\nu:B\sstac{L} H\to A$, $b\sstac{L} h\mapsto bj(h),$ which is
right $\hH$-colinear by the right colinearity of a cleaving map
and the left $B$-linearity of the 
{\color{blue}
$\hH_R$- and $\hH_L$-coactions}. 
The map $\nu$ is well
defined in view of (\ref{eq:js}). For all $b\sstac{L} h\in
B\sstac{L} H$,
\begin{eqnarray*}
\kappa\big(\nu(b\stac{L} h)\big) &=&
b{j(h)_{[0]}}^{[0]}\conv{j}\left( {j(h)_{[0]}}\ti\right)\stac{L}
        {j(h)_{[1]}}
= bj({h\di}\ui) \conv{j}( {h\di}\uii)\stac{L} {h}\dii\\
&=& b\eta_L\big( \pi_L(h\di)\big)\stac{L} h\dii = b\stac{L}
s_L\big(\pi_L(h\di)\big)h\dii=b\stac{L} h,
\end{eqnarray*}
where the second equality follows by the $\hH$-colinearity of $j$,
and the third one by (\ref{eq:rinv}). On the other hand,
(\ref{eq:linv}) and the counitality of $\roA$ imply, for all $a\in
A$,
$$
\nu\big( \kappa(a)\big)=a^{[0]}\conv{j}({a^{[1]}}\di)j({a^{[1]}}\dii)
=a^{[0]}\eta_R\big(\pi_R(a^{[1]})\big)=a.
$$
This means that $\nu$ is the left $B$-linear right $\hH$-colinear inverse of
$\kappa$, hence $\kappa$ is the required isomorphism.

(2) $\Ra$ (1) Suppose that the canonical map \eqref{eq:Rcan} is
bijective and write $\tau={\rm can}_R^{-1}(1_A\stac{R} -):H\to
A\sstac{B} A$ for the translation map. In explicit calculations we
use a Sweedler type notation, for all $h\in H$,
$\tau(h)=h^{\{1\}}\sstac{B} h^{\{2\}}$ (summation understood). Let
 $\kappa: A\to B\sstac{L} H$ be an isomorphism of left
$B$-modules and of right $\hH$-comodules and define maps $H\to A$
$$
j:= \kappa^{-1}(1_B\stac{L} -), \qquad
\tilde{j} := [A\stac{B} (B\stac{L}\pi_L)\circ\kappa]\circ
  \tau.
$$
We claim that $j$ is a cleaving map and $\tilde{j}$ is its
convolution inverse. First note that since $\kappa^{-1}$ is left
$B$-linear, it is in particular left $L$-linear, hence also $j$ is
left $L$-linear. Since $\kappa^{-1}$ is also right $\hH$-colinear,
so is $j$. Furthermore, the canonical map is left $A$-linear,
hence also left $R$-linear. Therefore, its inverse is left
$R$-linear, implying that, for all $h\in H$ and $r\in R$,
$\tau(h\, t_R(r))=\eta_R(r) h^{\{1\}}\sstac{B} h^{\{2\}}$. With
this property of the translation map at hand, one immediately
finds that, for all
 $h\in H$ and $r\in R$, ${\tilde j}(h t_R(r))
=\eta_R(r) {\tilde j}(h)$. On the other hand, by
\eqref{eq:robilin}, for all $a,a'\in A$ and $r\in R$, $ {\rm
can}_R(a\sstac{B} \eta_R(r) a')=aa^{\prime[0]}\sstac{R} s_R(r)
a^{\prime [1]}. $ This implies that $\tau(s_R(r)\,
  h)=h^{\{1\}}\sstac{B} \eta_R(r)h^{\{2\}}$. Thus, in view of the Hopf
  algebroid axiom $t_L = s_R\circ\pi_R\circ t_L$, one finds, for all $h\in H$
  and $l\in L$,
$$
{\tilde j}(t_L(l)\, h)
= [A\stac{B} (B\stac{L}\pi_L)\circ\kappa]\big( \tau(t_L(l) h)\big)
= [A\stac{B} (B\stac{L}\pi_L)\circ \kappa](h^{\{1\}}\stac{B}
\eta_R(\pi_R( t_L(l)))  h^{\{2\}}).
$$
  Since $\kappa$ is right
  $\hH_R$-colinear, it is in particular left $R$-linear, where the left
  $R$-module structure of $B\sstac{L} H$ is
  given by $r\cdot (b\sstac{L} h)=b\sstac{L} s_R(r) h$,
  (cf.\ \eqref{eq:RacM}). By
the right $L$-linearity of $\pi_L$ and the axiom
$t_L = s_R\circ\pi_R\circ t_L$, one therefore concludes that
${\tilde j}(t_L(l)\, h) = {\tilde j}(h) \eta_L(l)$, as required.
This proves that ${\tilde j}$ satisfies (\ref{eq:tjt}).
It remains to check (\ref{eq:rinv}) and (\ref{eq:linv}):
\begin{eqnarray*}
\mu_A\circ (j\stac{R} {\tilde j})\circ \gamma_R
&=& \mu_A\circ \{j\stac{R} [A\stac{B} (B\stac{L} \pi_L)\circ \kappa]\circ
\tau\}\circ
\gamma_R\\
&=& [A\stac{B} (B\stac{L} \pi_L)\circ\kappa]\circ {\rm can}_R^{-1}\circ
(j\stac{R}
H)\circ \gamma_R \\
&=& [A\stac{B} (B\stac{L} \pi_L)\circ\kappa]\circ {\rm can}_R^{-1}\circ
\roA\circ j\\
&=& [A\stac{B} (B\stac{L} \pi_L)\circ\kappa]\circ (1_A\stac{B} j(-))
=\eta_L\circ\pi_L,
\end{eqnarray*}
where the second equality follows by the left $A$-linearity of the canonical
map ${\rm can}_R$, hence of ${\rm can}_R^{-1}$, the third one by the right
$\hH_R$-colinearity of $j$ and the fourth one by the explicit form
\eqref{eq:Rcan} of ${\rm can}_R$. Furthermore,
\begin{eqnarray*}
\mu_A\circ ({\tilde j}\stac{L} {j})\circ \gamma_L
&=&  \mu_A\circ \{[A\stac{B} (B\stac{L} \pi_L)\circ \kappa]\circ \tau\stac{L}
j\}
\circ \gamma_L\\
&=&  \mu_A\circ[A\stac{B} (B\stac{L} \pi_L)\circ \kappa\stac{L} j]\circ
(A\stac{B} \lambda^A)\circ \tau\\
&=&  \mu_A\circ (A\stac{B} B\stac{L} \pi_L\stac{L} j)\circ(A\stac{B}
B\stac{L}\gamma_L)\circ (A\stac{B} \kappa)\circ \tau\\
&=&\mu_A\circ [A\stac{B} (B\stac{L} j)\circ \kappa]\circ \tau
=\mu_A\circ \tau=\eta_R\circ \pi_R,
\end{eqnarray*}
where the second equality follows by the $\hH_L$-colinearity of
$\tau$, the third one by the $\hH_L$-colinearity of $\kappa$, the
penultimate one by the left $B$-linearity of $\kappa$ and the last
one by $(A\stac{R}\pi_R)\circ {\rm can}_R=\mu_A$ and the
definition of the translation map $\tau$.
\end{proof}

By Remark \ref{rem:lefthand}, the following `left handed version' of
Theorem \ref{thm.main.cleft} (1)~$\Ra$~(2)(b) can be formulated.
\begin{corollary}\label{cor:kappaL}
Let $\hH$ be a Hopf algebroid with a bijective antipode and $B\subseteq A$ an
$\hH$-cleft extension with a cleaving map $j$. Then the  right $B$-linear
left $\hH$-colinear map
\begin{equation}\label{eq:kappaL}
A\to H\stac{L} B, \qquad
a\mapsto S^{-1}(a\bi)\di\stac{L} j\big(S^{-1}(a\bi)\dii\big) a\bz
\equiv S^{-1}(a\ti)\stac{L}j\big(S^{-1}({a\tz}\bi)\big){a\tz}\bz ,
\end{equation}
is an isomorphism.
\end{corollary}
The following is an immediate consequence of Theorem \ref{thm.main.cleft}.
\begin{corollary}\label{cor:ff}
Let $\hH$ be a Hopf algebroid and $B\subseteq A$ an $\hH$-cleft extension.
If $H$ is a
projective left $L$-module, then $A$ is a faithfully flat left $B$-module.
\end{corollary}
\begin{proof}
By Theorem~\ref{thm.main.cleft}~(1)~$\Ra$~(2)(b), $A\simeq
B\sstac{L} H$ as left $B$-modules. Since $H$ is projective as a
left $L$-module, $A$ is projective as a left $B$-module. Together
with Lemma \ref{lem:dirsum} this implies the claim.
\end{proof}

If the antipode of a Hopf algebroid $\hH$ is bijective then, by
\cite[Lemma 3.3]{Bohm:gal}, an extension $B\subseteq A$ is a right
$\hH_R$-Galois extension if and only if it is a right $\hH_L$-Galois
extension. By \cite[Lemma 4.1]{Bohm:gal}, this is further
equivalent to the left $\hH_R$-Galois and also to the
left $\hH_L$-Galois property of the extension. Hence repeating the
steps in the proof of \cite[Proposition 4.1]{BohBrz:str}, we conclude that
Lemma \ref{lem:dirsum}, Theorem \ref{thm.main.cleft} and Remark
\ref{rem:lefthand} imply the following
\begin{corollary}\label{cor:rel.inj}
Let $\hH$ be a Hopf algebroid with a bijective antipode
and $B\subseteq A$ an $\hH$-cleft extension.
Then $A$ is an $R$-relative injective right and left $\hH_R$-comodule, and
an $L$-relative injective left and $right \hH_L$-comodule.
\end{corollary}

\section{Crossed products with Hopf algebroids}\label{sec:crossp}

One of the main results in the theory of cleft extensions of Hopf algebras is
the
equivalent characterisation of such extensions as crossed product algebras
with an invertible cocycle (cf.\ \cite[Theorem~11]{DoiTak:cle}
\cite[Theorem~1.18]{BlaMon:cro}). The aim of this section is to derive such
a characterisation for Hopf algebroid cleft extensions. First we need to
develop
a suitable theory of crossed products, generalising that of \cite{DoiTak:cle}
and
\cite{BlaCoh:cro}. We start by extending the notion of a {\em measuring}
\cite[p.\ 139]{Swe:Hop}.
\begin{definition}\label{def:meas}
Let ${\mathcal L}=(H,L,s,t,\gamma,\pi)$ be a left bialgebroid and
$B$ an $L$-ring with unit map $\iota:L\to B$. ${\mathcal L}$ {\em
measures} $B$ if there exits a $k$-linear map, called a {\em
measuring}, $H\sstac{k} B\to B$,  $h\sstac{k} b\mapsto h\cdot b$
such that, for all $h\in H$, $l\in L$, $b,b'\in B$,
\begin{blist}
\item $h\cdot 1_B=\iota\big(\pi(h)\big)$;
\item $\big(t(l) h \big)\cdot b=(h\cdot b)\iota(l)$ and
$\big(s(l) h \big)\cdot b=\iota(l)(h\cdot b)$;
\item $h\cdot (b\,b')=(h\di \cdot b)(h\dii\cdot b')$.
\end{blist}
\end{definition}
Note that condition (b) means simply that a measuring is an $L$-$L$ bimodule
map, where $H$ is viewed as an $L$-$L$ bimodule via the left multiplication
by $s$ and $t$.
A left ${\mathcal L}$-module algebra $B$ is measured by
${\mathcal L}$ with a measuring provided by the left $H$-multiplication in $B$.
\begin{definition}\label{def:cocycle}
Let ${\mathcal L}=(H,L,s,t,\gamma,\pi)$ be a left bialgebroid and
$\iota:L\to B$ an $L$-ring, measured by ${\mathcal L}$. A
$B$-valued {\em 2-cocycle} $\sigma$ on ${\mathcal L}$ is a
$k$-linear map $H\sstac{L^{op}} H\to B$ (where the right and left
$L^{op}$-module structures on $H$ are given by right and left
multiplication by $t(l)$, respectively) satisfying
\renewcommand\descriptionlabel[1]
{\hspace{\labelsep}\textrm{#1}}
\begin{description}
\item[(a)] $\sigma(s(l)h,k)=\iota(l) \sigma(h,k)$ and
$\sigma(t(l)h,k)=\sigma(h,k)\iota(l)$;
\item[(b)]
$(h\di\cdot\iota(l))\sigma(h\dii,k)=\sigma(h,s(l)k)$;
\item[(c)] $\sigma(1,h)=\iota\big(\pi(h)\big)=\sigma(h,1)$;
\item[(d)] $[h\di\cdot \sigma(k\di,m\di)]\,\sigma(h\dii,k\dii m\dii)=
\sigma(h\di,k\di)\, \sigma(h\dii k\dii, m),
$
\end{description}
for all $h,k,m\in H$, $l\in L$.

An ${\mathcal L}$-measured $L$-ring $B$ is called a {\em $\sigma$-twisted left
${\mathcal L}$-module} if a $2$-cocycle $\sigma$ satisfies
\begin{description}
\item[(e)] $1_H\cdot b=b$,
\item[(f)] $[h\di \cdot(k\di\cdot b)]\,\sigma(h\dii,k\dii)=
\sigma(h\di,k\di)\,(h\dii k\dii\cdot b)$,
\end{description}
for all $h,k\in H$, $b\in B$.
\end{definition}
Conditions (c) in Definition~\ref{def:cocycle}
determine the normalisation of  $\sigma$ and  (d)  is a cocycle
condition. These
have the same form as corresponding conditions in the bialgebra case.
Conditions (a) determine the module map properties of $\sigma$ while (b)
ensures that $\sigma$ is properly $L$-balanced; both are needed for (d) (and
(f)) to make
sense. Condition (e) sates that a measuring is a {\em weak action} (cf.\
\cite[Definition~1.1]{BlaCoh:cro}).

 Similarly to the bialgebra case,
 the map $\sigma(h,h'):= \iota\big(\pi(h\,h')\big)$ is a (trivial) cocycle
  for an ${\mathcal L}$-measured $L$-ring $B$ with unit
$\iota$, provided that the measuring restricts to the action on $L$,
 $h\cdot \iota(l)=\iota\big(\pi(hs(l))\big)$, for $h\in H$ and $l\in L$.
A twisted left ${\mathcal L}$-module  corresponding to this trivial cocycle
 $\sigma$
is simply a left ${\mathcal L}$-module algebra.

\begin{proposition}\label{lem:excp}
Let ${\mathcal L}=(H,L,s,t,\gamma,\pi)$ be a left bialgebroid and
$\iota:L\to B$ an $L$-ring, measured by ${\mathcal L}$. Let
$\sigma:H\sstac{L^{op}} H\to B$ be a map, satisfying properties
(a) and (b) in Definition~\ref{def:cocycle}. Consider the
$k$-module $B\sstac{L} H$, where the left $L$-module structure on
$H$ is given by multiplication by $s(l)$ on the left. $B\sstac{L}
H$ is  an associative algebra with unit $1_B\sstac{L} 1_H$ and
product
\begin{equation}\label{sigma.prod}
(B\stac{L} H)\stac{k} (B\stac{L} H)\to (B\stac{L} H),\quad
(b\stac{L} h)\stac{k}(b'\stac{L} h')\mapsto
b(h\di\!\cdot\! b')\sigma(h\dii, h'\di)\stac{L} h_{(3)} h'\dii ,
\end{equation}
if and only if $\sigma$ is a cocycle and
$B$ is a $\sigma$-twisted ${\mathcal L}$-module.
The resulting associative algebra is called a {\em crossed product} of $B$
with ${\mathcal L}$ and is denoted by $B\cp _{\sigma} \cL$.
\end{proposition}
Note that the smash product algebra in Example \ref{ex:smash} is a crossed
product with a trivial cocycle.

\begin{proof} The element $1_B\cp  1_H$ is a left unit if and only if
\begin{equation}\label{eq:LU}
b\cp h=(1_H\cdot b)\sigma(1_H,h\di)\cp  h\dii, \qquad \textrm{for all }\ b\cp
h\in B\stac{L} H.
\end{equation}
If $\sigma(1_H,h)=\iota\big(\pi(h)\big)$ and $1_H\cdot b=b$, then
\eqref{eq:LU} obviously holds. On the other hand, applying
$B\sstac{L}\pi$ to \eqref{eq:LU} we arrive at the identity
\begin{equation}\label{eq:appi}
b\iota\big(\pi(h)\big)=(1_H\cdot b)\, \sigma(1_H,h), \qquad \textrm{for all
}\ b\in B,\ h\in H.
\end{equation}
Setting $b=1_B$ in \eqref{eq:appi} we obtain
$\sigma(1_H,h)=\iota\big(\pi(h)\big)$,  and setting
$h=1_H$ we get $1_H \cdot b=b$.
Analogously, the condition that $1_B\cp  1_H$ is a right unit is
equivalent to the condition $\sigma(h,1_H)=\iota\big(\pi(h)\big)$, for
all $h\in H$.

The associative law for product \eqref{sigma.prod}
 reads, for all $h,k,m\in H$, $a,b,c\in B$,
\begin{eqnarray}
&&a(h\di\cdot b)\sigma(h\dii,k\di)(h_{(3)}k\dii\cdot c)
\sigma(h_{(4)}k_{(3)},m\di)\cp  h_{(5)}k_{(4)}m\dii=\nonumber\\
&&a(h\di\cdot b)[h\dii\cdot(k\di\cdot c)]
[h_{(3)}\cdot \sigma(k\dii,m\di)] \sigma(h_{(4)},k_{(3)}m\dii)\cp
h_{(5)} k_{(4)} m_{(3)}.\label{eq:ass}
\end{eqnarray}
If $\sigma$ is a cocycle and $B$ is a $\sigma$-twisted module, then
\eqref{eq:ass} obviously holds. Note that,
for all $h,k\in H$,
\begin{equation}
\sigma(h\di,k\di)\iota\big(\pi(h\dii k\dii)\big)
=\sigma(h,k).
\label{eq:pisi}
\end{equation}
Applying $B\sstac{L}\pi$ to \eqref{eq:ass}, using \eqref{eq:pisi}
 and setting $a=1_B=b$, we arrive
at
\begin{equation}\label{eq:a1b1}
\sigma(h\di,k\di)(h\dii k\dii\!\cdot\! c)\sigma(h_{(3)}k_{(3)},\! m) \!=\!
[h\di\!\cdot\! (k\di\!\cdot\! c)\sigma(k\dii,\! m\di)]
\sigma(h\dii,k_{(3)}m\dii).
\end{equation}
Setting $c=1_B$ in
\eqref{eq:a1b1} we derive the cocycle condition
Definition~\ref{def:cocycle}~(d), while
setting $m=1_H$ in \eqref{eq:a1b1} we obtain
 Definition~\ref{def:cocycle}~(f).
\end{proof}

\begin{theorem}\label{thm.cross}
Let ${\mathcal L}=(H,L,s,t,\gamma,\pi)$ be a left bialgebroid and
$\iota:L\to B$ an $L$-ring. View $A=B\sstac{L}H$ as a left
$B$-module and a right ${\mathcal L}$-comodule in canonical ways
(i.e.\ the left $B$-multiplication is given by product in $B$ and
the right ${\mathcal L}$-coaction is $B\sstac{L}\gamma$, with the $L$-actions
on $H$ given by the left multiplication by $s$ and $t$). Then $A$
is a right ${\mathcal L}$-comodule algebra with unit
$1_B\sstac{L}1_H$ and a left $B$-linear product if and only if $A$
is a crossed product algebra as in Proposition~\ref{lem:excp}.
\end{theorem}
\begin{proof}
The definition of  the product in $B\cp _{\sigma} \cL$ immediately
implies that $B\cp _{\sigma} \cL$ is a right ${\mathcal
L}$-comodule algebra with a left $B$-linear multiplication.
Conversely, suppose that $A$ has the required ${\mathcal
L}$-comodule algebra structure. Then, in particular, $A$ is an
$L^{op}$-ring via $ L^{op}\to A$, $l\mapsto 1_B\sstac{L}t(l)$. We
use the hom-tensor relation
\begin{equation}\label{hom-ten}
\lRhom B{\mathcal L}{(B\stac{L}H)\stac{L^{op}}
(B\stac{L}H)}{ B\stac{L}H}\simeq \lrhom LL
{H\stac{L^{op}}(B\stac{L}H)}{ B}
\end{equation}
and the ${\mathcal L}$-colinearity of the product in $A$, to  view
the multiplication in $A$ as an $L$-$L$ bilinear map
$H\sstac{L^{op}}(B\sstac{L}H) \to B$. For any $b\in B$ and $h\in
H$, define
\begin{equation}\label{mea.a}
h\cdot b := (B\stac{L}\pi)\big( (1_B\stac{L}h)(b\stac{L}1_H)\big).
\end{equation}
By \eqref{hom-ten}, the above definition implies that, conversely,
\begin{equation}\label{a.mea}
(1_B\stac{L}h)(b\stac{L}1_H) = h\sw 1\!\cdot\! b\,\stac{L}\, h\sw 2.
\end{equation}
Now, the assumption that $1_B\sstac{L} 1_H$ is the unit in $A$
implies condition (a) in Definition~\ref{def:meas}. The conditions
(b) follow by the right $L$-linearity and  the left $B$-linearity
of the product respectively (remember that every right ${\mathcal
L}$-comodule map is necessarily right $L$-linear). The condition
(c) follows by the associativity of the product. Thus $B$ is
measured by $\cL$ with measuring \eqref{mea.a}.

Next, for all $h,h'\in H$, define
\begin{equation}\label{cyc.a}
\sigma(h,h') := (B\stac{L}\pi)\big((1_B\stac{L}h)(1_B\stac{L}h')\big).
\end{equation}
Then, by \eqref{hom-ten},
\begin{equation}\label{a.cyc}
(1_B\stac{L}h)(1_B\stac{L}h') = \sigma(h\sw 1, h'\sw 1)\stac{L}h\sw 2
h'\sw 2.
\end{equation}
Since $A$ is an $L^{op}$-ring, \eqref{cyc.a} defines a $k$-linear
map $\sigma: H\sstac{L^{op}}H\to B$. The conditions (a) in
Definition~\ref{def:cocycle}
 follow by the left $B$-linearity
and  the right $L$-linearity of the product respectively. To check condition
(b), take any $h,k\in H$ and $l\in L$, and compute
\begin{eqnarray*}
(h\di\!\cdot\!\iota(l))\,\sigma(h\dii,k) \!\!\!&=& \!\!\!
(B\stac{L}\pi)[(1_B\stac{L}h)(\iota(l)\stac{L}1_H)(1_B\stac{L}k)]\\
&=& \!\!\! (B\stac{L}\pi)[(1_B\stac{L}h)(1_B\stac{L}s(l)k)] =
\sigma(h,s(l)k),
\end{eqnarray*}
where the first and last equalities follow by the definitions  of
the measuring
and $\sigma$ and equations \eqref{a.mea}, \eqref{a.cyc},
and the left $B$-linearity
of the product. Finally, for all $b,b'\in B$, $h,h'\in H$,
\begin{eqnarray*}
(b\stac{L} h)(b'\stac{L} h') &=& b[(1_B\stac{L} h)(b'\stac{L}
    1_H)(1_B\stac{L}h')]\\
&=& b(h\sw 1\cdot b')[(1_B\stac{L} h\sw 2)(1_B\stac{L}h')] =
b\,(h\di\cdot b')\,\sigma(h\dii,h'\di)\stac{L} h_{(3)} h'\dii,
\end{eqnarray*}
where we have used the left $B$-linearity of the product and
\eqref{a.mea} and  \eqref{a.cyc}.
Proposition~\ref{lem:excp} yields the assertion.
\end{proof}

\begin{corollary}\label{cor.gauge}
Given a crossed product $B\cp_\sigma\cL$ and a convolution invertible map
$\chi\in \lrhom LLHB$ such that $\chi(1_H)=1_B$,
define, for all $h,k\in H$ and $b\in B$,
\begin{equation}\label{gauge.mea}
h\cdot^\chi b := \chi(h_{(1)})(h_{(2)}\cdot b)\conv{\chi}(h_{(3)}),
\end{equation}
\begin{equation}\label{gauge.coc}
\sigma^\chi(h,k) := \chi(h_{(1)})(h_{(2)}\cdot \chi(k_{(1)}))\sigma(h_{(3)},
k_{(2)})\conv{\chi}(h_{(4)}k_{(3)}).
\end{equation}
Then $B$ is a $\sigma^\chi$-twisted $\cL$-module with measuring
\eqref{gauge.mea}. The corresponding crossed product
$B\cp_{\sigma^\chi}\cL$ is called a {\em
gauge transform of $B\cp_\sigma\cL$}.
\end{corollary}
\begin{proof}
Any convolution invertible map $\chi\in \lrhom LLHB$ defines a
left $B$-module right ${\mathcal L}$-comodule automorphism $\Phi$
of $B\sstac{L}H$, by
\begin{equation}\label{eq.phi.gauge}
\Phi(b\stac{L}h) = b{\chi}(h\sw 1)\stac{L}h\sw 2, \qquad
\Phi^{-1}(b\stac{L}h) = b\conv{\chi}(h\sw 1)\stac{L}h\sw 2.
\end{equation}
If $\chi(1_H) =1_B$, then $\Phi(1_H\sstac{L}1_B) =
1_H\sstac{L}1_B$. We can use this isomorphism to induce a new
right ${\mathcal L}$-comodule algebra structure on $B\sstac{L} H$
(with unit $1_H\sstac{L}1_B$) from that of $B\cp_\sigma\cL$. In
view of Theorem~\ref{thm.cross}, this necessarily is a crossed
product with the measuring and cocycle given by equations
\eqref{mea.a} and \eqref{cyc.a}, i.e., for all $b\in B$ and
$h,k\in H$,
$$
h\cdot^\chi b = (B\stac{L}\pi)\big(\Phi^{-1}
\big(\Phi(1_B\stac{L}h)\Phi(b\stac{L}1_H)\big)\big),
$$
$$
\sigma^\chi(h,k) = (B\stac{L}\pi)\big(\Phi^{-1}
\big(\Phi(1_B\stac{L}h)\Phi(1_B\stac{L}k)\big)\big),
$$
where the product is computed in $B\cp_\sigma\cL$. One easily checks that
these have the form stated in equations \eqref{gauge.mea} and
\eqref{gauge.coc}.
\end{proof}

\begin{definition}\label{def.equiv}
Let ${\mathcal L}=(H,L,s,t,\gamma,\pi)$ be a left bialgebroid and $B$
an
$L$-ring. Crossed products $B\cp_\sigma\cL$ and $B\cp_{\bar{\sigma}}\cL$
are said to be {\em equivalent} if there exists a left $B$-module
isomorphism of right ${\mathcal L}$-comodule algebras
$B\cp_{\bar{\sigma}}\cL \to B\cp_\sigma\cL$.
\end{definition}
\begin{theorem}\label{thm.equiv}
Let ${\mathcal L}=(H,L,s,t,\gamma,\pi)$ be a left bialgebroid and $B$
an
$L$-ring. Crossed products $B\cp_\sigma\cL$ and $B\cp_{\bar{\sigma}}\cL$
are equivalent if and only if $B\cp_{\bar{\sigma}}\cL$ is a gauge transform
of $B\cp_\sigma\cL$.
\end{theorem}
\begin{proof}
In view of the hom-tensor relation $\lRhom B{\mathcal
L}{B\sstac{L}H}{B\sstac{L}H}\simeq \lrhom LLHB$, there is a
bijective correspondence between left $B$-module right ${\mathcal
L}$-comodule isomorphisms $\Phi: B\cp_{\bar{\sigma}}\cL\to
B\cp_\sigma\cL$ and convolution invertible maps $\chi\in \lrhom
LLHB$. This correspondence is given by equation
\eqref{eq.phi.gauge} in one direction and by $\chi(h) =
(B\sstac{L}\pi)(\Phi(1_B\sstac{L}h))$ in the other. If $\Phi$ is
also an algebra map, then $\chi(1_H) = 1_B$ and, following the
same line of argument as in the proof of
Corollary~\ref{cor.gauge}, one finds that the measuring
corresponding to $\bar{\sigma}$ is given by $h\cdot^\chi b$ and
that $\bar{\sigma} = \sigma^\chi$. Conversely, given $\chi$ and
corresponding (by equations \eqref{eq.phi.gauge}) isomorphism
$\Phi: B\cp_{{\sigma^\chi}}\cL\to B\cp_\sigma\cL$, one can
compute, for all $b,b'\in B$, $h,h'\in H$,
\begin{eqnarray*}
\Phi((b\cp_{{\sigma^\chi}} h)
(b'\!\!\!\!\!\!\!\!\!\!\!\!\!\!&&\!\!\!\!\!\!\!\cp_{{\sigma^\chi}}h')) =
b(h\di\!\cdot\!^\chi b')\sigma^\chi(h\dii, h'\di)\chi(h_{(3)} h'\dii)
\cp_{{\sigma}} h_{(4)} h'\sw 3\\
&=&\!\!\!
 b\chi(h\sw 1)(h\sw 2\!\cdot\! b')(h\sw 3\!\cdot\!
  \chi(h'\sw 1))
\sigma(h\sw 4,h'\sw 2)
\cp_{{\sigma}}
h\sw 5h'\sw 3\\
&=&\!\!\!
 b\chi(h\sw 1)(h\sw 2\!\cdot\! (b'
  \chi(h'\sw 1))
\sigma(h\sw 3,h'\sw 2)
 \cp_{{\sigma}}
h\sw 4h'\sw 3 = \Phi(b\cp_{{\sigma^\chi}} h)\Phi(b'\cp_{{\sigma^\chi}}h'),
\end{eqnarray*}
where the second equality follows by the fact that $\conv{\chi}$ is the
convolution inverse of $\chi$ and the counit axiom, and the third equality
follows by property (c) in
Definition~\ref{def:meas}. This proves that $\Phi$ is an algebra map, hence the
crossed product algebras $B\cp_{{\sigma^\chi}}\cL$ and $
B\cp_\sigma\cL$ are mutually equivalent.
\end{proof}

Next we establish what is meant by an invertible cocycle in this
generalised context.
\begin{definition}\label{def:cocinv}
Let ${\mathcal L}=(H,L,s,t,\gamma,\pi)$ be a left bialgebroid and
$\iota:L\to B$ an $L$-ring, measured by ${\mathcal L}$. A
$B$-valued 2-cocycle $\sigma$ on ${\mathcal L}$ is {\em
invertible} if there exists a $k$-linear map ${\tilde\sigma}:
H\sstac{L} H\to B$ (where the right and left $L$-module structures
on $H$ are given by right and left multiplication by $s(l)$,
respectively) satisfying
\begin{blist}
\item ${\tilde \sigma}(s(l)h,k)=\iota(l) {\tilde \sigma}(h,k)$ and
${\tilde \sigma}(t(l)h,k)={\tilde\sigma}(h,k)\iota(l)$;
\item
${\tilde \sigma}(h\di,k)(h\dii \cdot \iota(l))={\tilde \sigma}(h, t(l)k)$;
\item $\sigma(h\di,k\di)\,{\tilde \sigma}(h\dii,k\dii)=h\cdot(k\cdot 1_B)$ and
${\tilde \sigma}(h\di,k\di)\,\sigma(h\dii,k\dii)=hk\cdot 1_B$,
\end{blist}
for all $h,k\in H$ and $l\in L$. A map $\tilde{\sigma}$ is called an {\em
  inverse}
of $\sigma$.
\end{definition}
Again, conditions (a) and (b) are needed so that the inverse property (c)
can be stated. In the case $\cL$ is a bialgebra over a ring $L=k$,
conditions (a) and (b) are satisfied automatically. The following
two lemmas explore the nature of
cocycles and their inverses.
\begin{lemma}\label{lem:tsinorm}
Let ${\mathcal L}=(H,L,s,t,\gamma,\pi)$ be a left bialgebroid and $B$ an
$L$-ring with unit $\iota:L\to B$. Assume that ${\mathcal L}$ measures
$B$ and $\sigma$ is an invertible $B$-valued 2-cocycle on ${\mathcal
L}$. Then an inverse
$\tilde{\sigma}$ of $\sigma$ is unique and normalised, i.e.,
for all $h\in H$,
$$
{\tilde \sigma}(1_H,h)=\iota\big(\pi(h)\big)={\tilde \sigma}(h,1_H).
$$
\end{lemma}
\begin{proof}
Note that, if $\tilde{\sigma}$ is an inverse of $\sigma$, then, for all
$h,k\in H$,
\begin{equation}
\iota\big(\pi(h\di k\di)\big) {\tilde \sigma}(h\dii,k\dii)
={\tilde \sigma}(h,k).
\label{eq:pitsi}
\end{equation}
Using this identity and Definition~\ref{def:cocinv}~(c), one finds that
\begin{equation}
{\tilde \sigma}(h,k) = {\tilde\sigma}(h\di,k\di)[h\dii\cdot (k\dii \cdot 1_B)].
\label{eq:twb}
\end{equation}
Now suppose that $\hat{\sigma}$ is another inverse of $\sigma$. Then replacing
the expression in square brackets in \eqref{eq:twb} by the first of equations
in
Definition~\ref{def:cocinv}~(c) for $\hat{\sigma}$, using the second of
equations Definition~\ref{def:cocinv}~(c) for $\tilde{\sigma}$, and finally
using
\eqref{eq:pitsi} for $\hat{\sigma}$, one finds that $\tilde{\sigma} =
\hat{\sigma}$.
Hence the inverse of a cocycle is unique.

Use \eqref{eq:pitsi}, Definition \ref{def:cocycle} and
Definition~\ref{def:cocinv}~(c) to compute, for all $h\in H$,
\begin{eqnarray*}
{\tilde \sigma}(1_H,h)&=&
\iota\big(\pi(h\di)\big){\tilde \sigma}(1_H, h\dii)
=\sigma(1_H,h\di){\tilde \sigma}(1_H, h\dii)\\
&=&[1_H\cdot(h\cdot 1_B)]\sigma(1_H,1_H)
=\sigma(1_H, s(\pi(h)))
=\iota\big(\pi(h)\big).
\end{eqnarray*}
The proof of the other identity is similar.
\end{proof}
\begin{lemma}\label{lem:aconsi}
Let ${\mathcal L}=(H,L,s,t,\gamma,\pi)$ be a left bialgebroid, $B$ an
$L$-ring, measured by ${\mathcal L}$, and $\sigma$ an
 invertible $B$-valued 2-cocycle on ${\mathcal L}$ with
 the inverse $\tilde{\sigma}$. For all $h,k,m\in H$,
\begin{blist}
\item
$h\cdot \sigma(k,m)=\sigma(h\di,k\di)\sigma(h\dii k\dii,m\di)
{\tilde \sigma}(h_{(3)},k_{(3)}m\dii)$,
\item
$
h\cdot {\tilde \sigma}(k,m)=\sigma(h\di,k\di m\di)
{\tilde \sigma}(h\dii k\dii,m\dii) {\tilde \sigma}(h_{(3)},k_{(3)}).
$
\end{blist}
\end{lemma}
\begin{proof}
(a)
Denote the unit map of the $L$-ring $B$ by $\iota:L\to B$.
In view of  \eqref{eq:pisi} and with the help of
properties (c) and (a) in Definition
\ref{def:meas} and  (c) in
Definition~\ref{def:cocinv}, we can
compute, for all $h,k,m\in H$,
\begin{eqnarray*}
h\!\cdot\! \sigma(k,m)
\!\!\!&=& \!\!\! h\!\cdot\! [\sigma(k\di,m\di)\iota\big(\pi(k\dii m\dii)\big)]
= [h\di \!\cdot\! \sigma(k\di,m\di)][h\dii \!\cdot\! \big(k\dii m\dii \!\cdot
1_B\big)]\\
&=& [h\di \!\cdot \sigma(k\di,m\di)]\sigma(h\dii,k\dii m\dii)
{\tilde \sigma}(h_{(3)}, k_{(3)} m_{(3)})\\
&=& \sigma(h\di,k\di)\sigma(h\dii k\dii,m\di)
{\tilde \sigma}(h_{(3)}, k_{(3)} m\dii),
\end{eqnarray*}
where the last equality follows by property (d) in Definition
\ref{def:cocycle}.

{{(b)}}
Use  part (a), \eqref{eq:pisi} and Definition~\ref{def:cocinv}~(c) to
find that, for all $h,k,m\in H$,
\begin{eqnarray}
[h\di\cdot \sigma(k\di,m\di)]
\!\!\!\!\!&\sigma&\!\!\!\!\!(h_{(2)}, k_{(2)} m_{(2)})
{\tilde\sigma}(h_{(3)} k_{(3)},m_{(3)})
{\tilde\sigma}(h_{(4)}, k_{(4)})\nonumber\\
&=&\sigma(h\di,k\di)
[h\dii k\dii \cdot(m\cdot 1_B)]
{\tilde\sigma}(h_{(3)}, k_{(3)}).\nonumber
\end{eqnarray}
By Definitions~\ref{def:meas} and~\ref{def:cocycle},
 $\sigma(h,s(\pi(m)))=h\cdot(m\cdot 1_B))$, for all $h,m\in H$. Hence
 conditions (d) in Definition \ref{def:cocycle} and (c) in Definition
 \ref{def:cocinv} allow us to develop the right
 hand
 side of the above equality further to arrive at the equation
 \begin{equation}
[h\di\cdot \sigma(k\di,m\di)]
\sigma (h_{(2)}, k_{(2)} m_{(2)})
{\tilde\sigma}(h_{(3)} k_{(3)},m_{(3)})
{\tilde\sigma}(h_{(4)}, k_{(4)})
=h\!\cdot\![k\!\cdot\!(m\!\cdot\! 1_B)].
\label{eq:auxi}
\end{equation}
Therefore
\begin{eqnarray*}
&\!\!\! h&\!\!\!\!\!\!\! \cdot {\tilde\sigma}(k,m)
=h\!\cdot\! \{{\tilde\sigma}(k\di,m\di)[k\dii\!\cdot\! (m\dii\!\cdot\!
  1_B)]\}\\
&=&[h\di\!\cdot\! {\tilde\sigma}(k\di,m\di)]
\{ h\dii\!\cdot\! [k\dii\!\cdot\! (m\dii\!\cdot\! 1_B)]\}\\
&=&[h\di\!\cdot {\tilde\sigma}(k\di,m\di)]
[h\dii\!\cdot \sigma(k\dii,m\dii)]
\sigma(h_{(3)}, k_{(3)} m_{(3)})
{\tilde\sigma}(h_{(4)} k_{(4)},m_{(4)})
{\tilde\sigma}(h_{(5)}, k_{(5)})\\
&=&[h\di\cdot(k\di m\di\cdot 1_B)]
\sigma(h_{(2)}, k_{(2)} m_{(2)})
{\tilde\sigma}(h_{(3)} k_{(3)},m_{(3)})
{\tilde\sigma}(h_{(4)}, k_{(4)})\\
&=&\sigma(h\di,k\di m\di){\tilde\sigma}(h_{(2)} k_{(2)},m_{(2)})
{\tilde\sigma}(h_{(3)}, k_{(3)}),
\end{eqnarray*}
where the first equality follows by \eqref{eq:twb}, the second one by
property (c) in Definition \ref{def:meas} and the third one by
\eqref{eq:auxi}.
The penultimate equality follows by property
(c) in Definition \ref{def:meas} and (c) in Definition
\ref{def:cocinv}. The last equality follows by
conditions (a) in Definition \ref{def:meas} and (b) in Definition
\ref{def:cocycle}.
\end{proof}

Take a Hopf algebroid $\hH$, an $\hH_L$-measured $L$-ring $B$ and
a cocycle $\sigma$. Then the crossed product $B\cp _\sigma \hH_L$
is an $R$-ring with unit map $\eta_R: r\mapsto 1_B\cp  s_R(r) $
and an $L$-ring with unit map $ \eta_L: l\mapsto 1_B\cp  s_L(l)
=\iota(l)\cp  1_H, $ where $\iota:L\to B$ denotes the unit map of
the $L$-ring $B$. It is also a right $\hH$-comodule with
$\hH_R$-coaction $B\sstac{L}\gamma_R$ and right
$\hH_L$-coaction $B\sstac{L}\gamma_L$. The 
{\color{blue}
$\hH_R$}-coinvariants are the
elements of the form $b\cp  1_H$ for $b\in B$, hence they form an
$L$-subring, isomorphic to $B$. Therefore, $B\subseteq B\cp
_\sigma \hH_L$  is a right $\hH$-extension, and it is natural to
ask whether this extension is cleft.
\begin{theorem}\label{thm.cross.cleft}
Let $\hH$ be a Hopf algebroid and $B\cp _\sigma \hH_L$ a crossed
product of an $\hH_L$-measured $L$-ring $B$. If the cocycle
$\sigma$ is  invertible,
then the extension $B\subseteq B\cp _\sigma \hH_L$ is $\hH$-cleft.
\end{theorem}
\begin{proof}
We claim that the right $\hH$-comodule
map  $j:H\to B\cp _\sigma \hH_L$, $h\mapsto
1_B\cp  h$ is a cleaving map with the convolution inverse
$$
\conv{j}(h)={\tilde\sigma}(S(h\di)\di,h\dii)\cp
S(h\di)\dii=
{\tilde\sigma}(S({h\uii}\di),{h\uii}\dii)\cp  S(h\ui).
$$
The two forms of $\conv{j}$ are equivalent by the
anti-comultiplicativity of $S$ and the left $\hH_R$-colinearity of
$\gamma_L$. Using the definitions of a cocycle and
its inverse, and in particular, the module and normalisation properties of
$\sigma$ and $\tilde{\sigma}$, one verifies that $j$ and
$\conv{j}$ have the required $L$-, $R$-module map 
properties
\eqref{eq:js} and \eqref{eq:tjt}.
Next, take  any $h\in H$ and compute
\begin{eqnarray*}
\conv{j}(h\di) j(h\dii)
&=& {\tilde \sigma}( S(h\di)\di,h\dii) \sigma(S(h\di)\dii,h_{(3)})\cp
S(h\di)_{(3)} h_{(4)}\\
&=&1_B\cp  S(h\di)h\dii=1_B\cp
s_R\big(\pi_R(h)\big)=\eta_R\big(\pi_R(h)\big),
\end{eqnarray*}
where the second equality follows by condition (c) in Definition
\ref{def:cocinv}, condition (a) in Definition \ref{def:meas} and the
counit property of $\pi_L$. The third equality follows by the antipode
axiom \eqref{ax:antip}. The proof of the  identity \eqref{eq:rinv} is slightly
more involved:
\begin{eqnarray*}
j(\!\!\!\!\!\!\!&h\ui&\!\!\!\!\!\!\!) \conv{j}(h\uii)\\
&=&[{h\ui}\di\cdot {\tilde \sigma}( S({h\uii}\di)\di,{h\uii}\dii)]
\sigma({h\ui}\dii, S({h\uii}\di)\dii)\cp
{h\ui}_{(3)}S({h\uii}\di)_{(3)}\\  
&=&{\sigma}({h\ui}\di,S({h\uii}\di)\di {h\uii}\dii)
{\tilde \sigma}({h\ui}\dii S({h\uii}\di)\dii, {h\uii}_{(3)})\cp
{h\ui}_{(3)} S({h\uii}\di)_{(3)}\\ 
&=&{\sigma}({h\ui}\di,
S({{{{h^{(3)}}\di}\uii}\di}) 
{{{{h^{(3)}}\di}\uii}\dii})
{\tilde \sigma}({h\ui}\dii 
S({{{h^{(3)}}\di}\ui}), 
{{h^{(3)}}\dii})\cp
{h\ui}_{(3)} S({h\uii})\\  
&=&{\tilde \sigma}\big( s_L\big(\pi_L({h\ui}\di 
s_R(\pi_R({{{h^{(3)}}\di}\uii})))\big)
{h\ui}\dii 
S({{{h^{(3)}}\di}\ui}),
{{h^{(3)}}\dii}\big)\cp {h\ui}_{(3)} S({h\uii})\\ 
&=&{\tilde \sigma}\big({h\ui}\di 
S({{{h^{(3)}}\di}\ui} s_R(\pi_R({{{h^{(3)}}\di}\uii}))),
{{h^{(3)}}\dii}\big)\cp {h\ui}_{(2)} S({h\uii})\\ 
&=&{\tilde \sigma}\big({h\ui}\di S({h^{(3)}}\di),{h^{(3)}}\dii\big)\cp
{h\ui}\dii S(h\uii)\\ 
&=&{\tilde \sigma}\big(
{ {h\ui}\di} 
S({{h\uii}\di}),
{ {h\uii}\dii}\big)\cp
{ {{h\ui}\dii}\ui} 
S({ {{h\ui}\dii}\uii})\\ 
&=&{\tilde \sigma}\big(t_L(\pi_L(
{ {h\ui}\dii}))
{ {h\ui}\di}
S(
{ {h\uii}\di}),
{ {h\uii}\dii}\big)\cp  1_H\\ 
&=&{\tilde \sigma}\big({h\di}\ui S({h\di}\uii), h\dii\big)\cp  1_H
={\tilde \sigma}(1_H,h)\cp  1_H=\eta_L\big(\pi_L(h)\big),
\end{eqnarray*}
where the second equality follows by Lemma \ref{lem:aconsi} (b),
condition (c) in Definition \ref{def:cocinv}, condition (a) in
Definition \ref{def:meas}, condition (a) in Definition
\ref{def:cocinv} and the counit property of $\pi_L$. The third
equality follows by the anti-comultiplicativity of $S$ and
\eqref{ax:cpcolin}. The fourth one follows by the antipode axiom
\eqref{ax:antip}, the fact that the domain of $\sigma$ is
$H\sstac{L^{op}} H$ (i.e. $\sigma(ht_L(l),k)=\sigma(h,t_L(l)k)$
for $h,k\in H$, $l\in L$), the normalisation of cocycles
(condition (c) in Definition \ref{def:cocycle}) and the left
$L$-linearity of ${\tilde \sigma}$ in the first argument
(condition (a) in Definition \ref{def:cocinv}). In the fifth step
the Hopf algebroid identity $\pi_L(hs_R(r))=\pi_L(hS(s_R(r)))$,
implying $s_L\big(\pi_L(h\di s_R(r))\big) h\dii = h S\big(
s_R(r)\big)$, for $h\in H$ and $r\in R$, has been used together
with the anti-multiplicativity of $S$. The sixth and seventh
equalities follow by the coassociativity and $\hH_L$-colinearity
of $\gamma_R$ and the counit property of $\pi_R$. The eighth
equality follows by the antipode axiom \eqref{ax:antip} and the
right $L$-linearity of ${\tilde \sigma}$. The ninth one follows by
axiom \eqref{ax:cpcolin} and the counit property of $\pi_L$. The
penultimate equality follows by axiom \eqref{ax:antip} and the
fact that the domain of ${\tilde \sigma}$ is $H\sstac{L} H$ (i.e.
${\tilde \sigma}(hs_L(l),k)={\tilde \sigma}(h,s_L(l)k)$ for
$h,k\in H$, $l\in L$). The last equality follows by Lemma
\ref{lem:tsinorm}.
Note that all expressions in the above computation are well defined. In order
to see that, recall that the restriction of $H\ot_R \pi_R$ to the Takeuchi
product $H\times_R H$ is right $L$-linear (making the expression in the fifth
line meaningful) and, similarly, the restriction of $H\ot_L \pi_L$ to
$H\times_L H$ is right $R$-linear (making the expression in the penultimate
line meaningful).
\end{proof}

The final  aim of this section is to prove that any cleft extension is
necessarily isomorphic
to a crossed product with an invertible cocycle.
\begin{theorem}\label{thm:crossed}
If $B\subseteq A$ is an $\hH$-cleft
extension, then there exists an invertible cocycle $\sigma$
and a left $B$-module right $\hH$-comodule algebra isomorphism
 $A\to B\cp _\sigma \hH_L$.
\end{theorem}
\begin{proof}
For an $\hH$-cleft extension $B\subseteq A$ the
cleaving map $j$ takes the unit element of $H$ to an invertible element of $B$
(with the inverse $\conv{j}(1_H)$). Thus, without the loss of generality,
 we can assume that a cleaving map $j$ is normalised, i.e.
$j(1_H)=1_A=\conv{j}(1_H)$. By Theorem~\ref{thm.main.cleft}, $A$
is isomorphic to $B\sstac{L} H$ as a left $B$-module and a right
$\hH$-comodule. We can use this isomorphism to induce a comodule
algebra structure on $B\sstac{L}H$. By Theorem~\ref{thm.cross},
the induced algebra structure is necessarily a crossed product
$B\cp _\sigma \hH_L$. In view of the definitions of the map
$\kappa$ and its inverse in the proof of
Theorem~\ref{thm.main.cleft}~(1)~$\Ra$~(2)(b), the measuring and
cocycle can be read off equations \eqref{mea.a} and \eqref{cyc.a},
respectively, and come out as
\begin{equation}\label{mea.coc}
h\cdot b = j(h\ui)\,b\,\conv{j}(h\uii),
\qquad \sigma(h,k) =  j(h\ui)\,j(k\ui) \,\conv{j}(h\uii k\uii).
\end{equation}
We only need to prove that the cocycle $\sigma$ is invertible.
Define
\begin{equation}\label{eq:jtotsi:a}
{\tilde \sigma}: H\stac{L} H\to B,\qquad
h\stac{L} k \mapsto j(h\ui k\ui)\,\conv{j}(k\uii)\,\conv{j}(h\uii).
\end{equation}
The map \eqref{eq:jtotsi:a} is well defined by \eqref{eq:js} and
\eqref{eq:tjt}, on one hand, and by \eqref{eq:tjt}, \eqref{eq:third}
and
the property that the range of the coproduct of a right bialgebroid is in
the Takeuchi product,
on the other hand. The proof
that the range of $\tilde{\sigma}$ is in $B$ and that
$\tilde{\sigma}$ is the inverse of the cocycle $\sigma$ are done
by a  routine calculation and are left to the reader.
\end{proof}

Combining Theorem~\ref{thm.equiv} with  Theorem~\ref{thm:crossed},
we can fully describe the relationship between different cleaving maps for
the same cleft extension.
\begin{corollary}\label{cor.gauge.cleft}
Let $\hH$ be a Hopf algebroid and $B\subseteq A$ an $\hH$-cleft
extension with a (non-necessarily unital)
cleaving map $j:H\to A$. Then a map $j':H\to A$ is a
cleaving map if and only if there exists a convolution invertible $L$-$L$
bilinear map
$\chi:H\to B$,
such that
\begin{equation}\label{eq:j'}
j'(h)=\chi(h\di)j(h\dii), \qquad {\textrm for\; all\ }h\in H.
\end{equation}
\end{corollary}
\begin{proof}
If $j$ is a cleaving map and $\chi \in {\rm Hom}_{L,L}(H,A)$ is convolution
invertible, then \eqref{eq:j'} obviously defines a cleaving map.
In order to prove the converse claim,
suppose first that both $j$ and $j'$ are normalised as in the proof of
Theorem~\ref{thm:crossed}. By
 Theorem~\ref{thm:crossed}, the crossed products corresponding to
$j$ and $j'$ are isomorphic to $A$ via left $B$-module right $\hH$-comodule
 algebra
maps, hence they are equivalent to each other. The isomorphism,
obtained from combining the maps $\nu$ (for $j'$) with $\kappa$
(for $j$) in the proof of
Theorem~\ref{thm.main.cleft}~(1)~$\Ra$~(2)(b), explicitly comes
out as $\Phi:b\sstac{L} h\mapsto bj'(h\sw 1\suc 1)\conv{j}(h\sw
1\suc 2)\sstac{L} h\sw 2$.
 Then, by Theorem~\ref{thm.equiv}, the existence of $\Phi$ is
 equivalent to the existence of a normalised convolution invertible map
$\chi\in \lrhom LLHB$, $\chi(h)= j'(h\ui)\conv{j}(h\uii)$.
Using the right $\hH_L$-colinearity of
$\gamma_R$ and the fact that $\conv{j}$ is a left convolution
inverse of $j$, one finds, for all $h\in H$, $\chi(h\di)j(h\dii)
= j'(h)$, i.e.\ equation~\eqref{eq:j'} holds. Allowing for $j$, $j'$ to be
non-unital
is equivalent to not requiring that $\chi$ be normalised.
\end{proof}

\section{The relative Chern-Galois character for $\hH$-cleft
extensions}
\label{sec.chern}

The aim of this section is to give a complete description of strong
connections in a cleft extension $B\subseteq A$ of a Hopf algebroid
$\hH = (\hH_L,\hH_R,S)$ (over rings $L$ and $R$)
and to find criteria for the existence and independence on the strong
connection
of the corresponding relative Chern-Galois characters introduced
and computed
in \cite{BohBrz:str}.

Begin with a right $\hH$-extension $B\subseteq A$ and suppose that
$T$ is a subalgebra of $B$. Then $A$ is called an {\em
$(\hH_R,T)$-projective left $B$-module} provided there exists a
left $B$-linear, right $\hH_R$-colinear section $\alpha_T$ of the
multiplication map $B\sstac{T}A\to A$. To consider the most
general case possible, we make no assumptions on a ring $T$ (but,
possibly, the most natural choice for $T$ is the base algebra
$L$).
\begin{lemma}\label{lem.rel.proj}
Let $\hH$ be a Hopf algebroid and $B\subseteq A$ an $\hH$-cleft
extension. Then
$A$ is an  $(\hH_R,L)$-projective left $B$-module.
\end{lemma}
\begin{proof}
The map $\tilde{\alpha}_L: B\sstac{L} H\to B\sstac{L}B\sstac{L}H$,
$b\sstac{L}h\mapsto b\sstac{L}1_B\sstac{L}h$ is a left $B$-linear
right $\hH$-colinear splitting of the product map
$b\sstac{L}b'\sstac{L} h\mapsto bb'\sstac{L}h$. By
Theorem~\ref{thm.main.cleft}, $A\simeq B\sstac{L} H$ as left
$B$-modules and right $\hH$-comodules, hence there is a
corresponding splitting $\alpha_L$ of the $B$-product map in $A$.
Explicitly,
$$
\alpha_L = (B\stac{L}\kappa^{-1})\circ\tilde{\alpha}_L\circ\kappa, \quad
a\ \mapsto
{a_{[0]}}^{[0]} \conv{j} ({a_{[0]}}^{[1]})\stac{L} j(a_{[1]})=
{a}^{[0]}\conv{j} ({a^{[1]}}_{(1)})\stac{L} j({a^{[1]}}_{(2)}),
$$
where $\kappa$ is the isomorphism \eqref{eq:kappa} in the proof of
Theorem~\ref{thm.main.cleft}
and $j$ is a cleaving map with the convolution inverse $\conv{j}$.
\end{proof}

Any right $\hH$-comodule algebra $A$ with 
{\color{blue}
$\hH_R$-coaction $\roA:a\mapsto a^{[0]}\ot_R a^{[1]}$ and 
$\hH_L$-coaction $\lambda^A:a\mapsto a_{[0]}\ot_L a_{[1]}$}
gives rise to an entwining map (over $R$) $\psi: H\sstac{R}A\to
A\sstac{R} H$, $h\sstac{R}a\mapsto a\su 0\sstac{R}ha\su 1$. The
map $\psi$ is bijective, provided the antipode $S$ is bijective
(cf.\ \cite[Lemma~4.1]{Bohm:gal}), and then the corresponding left
$\hH_R$-coaction on $A$ is
\begin{equation}\label{eq:Aro}
\Aro: A\to H\stac{R} A,\qquad a\mapsto S^{-1}(a\bi)\stac{R} a\bz
\end{equation}
(compare with \eqref{eq:L_R_coac}). Thus, following
\cite[Definition~3.4]{BohBrz:str}, if $B\subseteq A$ is a right
$\hH$-extension and $T$ is a subalgebra of $B$, then a left and
right $\hH_R$-comodule map $\ell_T: H\to A\sstac{T}A$ is a {\em
strong $T$-connection} provided that ${\widetilde {\rm
can}}_T\big(\ell_T(h)\big)=1_A\sstac{R} h$, for all $h\in H$,
where the map
\begin{equation}\label{eq:tcan}
{\widetilde {\rm can}}_T:A\stac{T} A\to A\stac{R} H,\qquad
a\stac{T} a'\mapsto aa^{\prime[0]}\stac{R} a^{\prime [1]},
\end{equation}
is well defined by the left $T$-linearity of $\roA$. The
$\hH_R$-coactions in $A\sstac{T}A$ are $A\sstac{T}\roA$ and
$\Aro\sstac{T}A$, with $\Aro$ given in \eqref{eq:Aro}.

The first observation is that a cleft extension comes equipped with a
strong $L$-connection.
\begin{theorem}\label{thm.str.con.can}
Let $\hH$ be a Hopf algebroid with a bijective antipode and $B\subseteq A$ an
$\hH$-cleft extension with a cleaving map $j$. Then the map
\begin{equation}\label{ell.can}
\ell_L: H\to A\stac{L} A,\qquad h\mapsto \conv{j}(h\di)\stac{L} j(h\dii),
\end{equation}
is a strong $L$-connection.
\end{theorem}
\begin{proof}
By Theorem~\ref{thm.main.cleft}, $B\subseteq A$ is a Galois $\hH_R$-extension,
which is $(\hH_R,L)$-projective by Lemma~\ref{lem.rel.proj}. Thus the existence
of a strong connection follows by \cite[Theorem~3.7]{BohBrz:str}. Using the
explicit forms of the inverse of the canonical $\hH_R$-Galois map in the
proof
of Theorem~\ref{thm.main.cleft} and of $\alpha_L$ in the proof
of Lemma~\ref{lem.rel.proj}, following the proof of
\cite[Theorem~3.7]{BohBrz:str} one arrives at the form of a strong
$L$-connection in \eqref{ell.can}. 
\end{proof}

{\color{blue}
The full classification
of strong $T$-connections in a cleft extension is described 
in forthcoming Theorem \ref{thm.class.con}. Its proof relies on two lemmata:

\begin{lemma}\label{lem:new_adj}
For any Hopf algebroid $\hH$, the following statements hold.

(1) The forgetful functor $\M^\hH \to \M_L$ possesses a
  right adjoint $(-)\ot_L H$.

(2) The forgetful functor ${}^\hH\M \to {}_L \M$ possesses a
  right adjoint $H \ot_L (-)$.
\end{lemma}

\begin{proof}
(1) The unit of the adjunction is given by the $\hH_L$-coaction $M \to M\ot_L
  H$, for any right $\hH$-comodule $M$. It is an $\hH$-comodule map by
  definition. The counit is given by $N\ot_L \pi_L:N\ot_L H \to N$, for any
  right $L$-module $N$. 
Part (2) is proven symmetrically.
\end{proof}

\begin{lemma}\label{lem:H-H_R}
Let $\hH$ be a Hopf algebroid with a bijective antipode and
$B\subseteq A$ an $\hH$-cleft extension. Then the obvious inclusion 
$\Hom^{\hH,\hH}(H, A\sstac T A)\hookrightarrow \Hom^{\hH_R,\hH_R}(H, A\sstac T
A)$ is an isomorphism.  
\end{lemma}

\begin{proof}
We need to show that any $f\in \Hom^{\hH_R,\hH_R}(H, A\sstac T A)$ is also
left and right $\hH_L$-colinear. Indeed, in terms of a cleaving map $j$ and
its convolution inverse $\conv{j}$, for any $h\in H$, 
\begin{eqnarray*}
f(h)&=&
f\big(h^{(1)}s_R(\pi_R(h^{(2)}))\big) =
f(h^{(1)}) \eta_R(\pi_R(h^{(2)})) =
f(h^{(1)}) \conv{j}({h^{(2)}}_{(1)}) j({h^{(2)}}_{(2)})\\
&=& f({h_{(1)}}^{(1)}) \conv{j}({h_{(1)}}^{(2)}) j(h_{(2)})=
f({h_{(1)}})^{[0]} \conv{j}(f({h_{(1)}})^{[1]}) j(h_{(2)}).
\end{eqnarray*}
The first equality follows by the counit property of $\pi_R$ and the second
one follows by the right $R$-linearity of $f$. In the third equality we used
that $\conv{j}$ is left convolution inverse of $j$, 
i.e. \eqref{eq:linv}. The penultimate equality is a
consequence of the right $\hH_L$-colinearity of $\gamma_R$ and the last
equality follows by the right $\hH_R$-colinearity of $f$. 
By Lemma \ref{lem:coinv}, the left $B$-linearity of the right $\hH_L$-coaction
on $A$ and the right $\hH_L$-colinearity of $j$, this proves
that $f$ is right $\hH_L$-colinear. Left  $\hH_L$-colinearity of $f$ is proven 
symmetrically. 
\end{proof}
}

\begin{theorem}\label{thm.class.con}
Let $\hH$ be a 
Hopf algebroid with a bijective antipode and
$B\subseteq A$ an $\hH$-cleft extension, and let $T$ be a
subalgebra of $B$. Write $\mu_B$ for the multiplication map
$B\sstac{T}B\to B$. Strong $T$-connections in $B\subseteq A$ are
in bijective correspondence with $L$-$L$ bilinear maps $f: H\to
B\sstac{T} B$ such that $\mu_B\circ f=\eta_L\circ \pi_L$.
\end{theorem}
\begin{proof}
Let $j$ be a cleaving map and $\conv{j}$ its convolution inverse.
By 
{\color{blue} Lemma \ref{lem:H-H_R}},
Theorem~\ref{thm.main.cleft} and
Corollary~\ref{cor:kappaL}\
{\color{blue} and finally Lemma \ref{lem:new_adj}},
there is a chain of isomorphisms
$$
\Bihom{\hH_R}{\hH_R} H {A\stac{T}A} 
{\color{blue}
\simeq \Bihom{\hH}{\hH} H {A\stac{T}A}}
\simeq \Bihom{\color{blue}\hH}{\color{blue} \hH} H 
{H\stac{L}B\stac{T}B\stac{L}H} \simeq \lrhom LLH{B\stac{T}B}.
$$
In view of the
explicit form of the isomorphism $\kappa: A\to B\sstac{L}H$ in
\eqref{eq:kappa} and its left-handed version \eqref{eq:kappaL}, we
thus obtain:
\begin{equation}\label{eq:f}
\Bihom{\hH_R}{\hH_R} H {A\stac{T}A}\ni\ell_T\mapsto [f:h\mapsto
  j(h\ui)\ell_T(h\uii) \conv{j}(h^{(3)})],
\end{equation}
and its inverse
\begin{equation}\label{eq:strcon}
\lrhom LLH{B\stac{T}B}\ni f\mapsto [\ell_T: h \mapsto \conv{j}(h\di) f(h\dii)
j(h_{(3)})].
\end{equation}
If $\ell_T$ is  a strong $T$-connection, then  $\mu_A\circ
\ell_T=\eta_R\circ \pi_R$, where $\mu_A:A\sstac{T}A\to A$ is the
product map in $A$. This implies that for the corresponding $f$ in
\eqref{eq:f}, $\mu_B\circ f=\eta_L\circ \pi_L$. Conversely,
suppose that $f$ has this property, write $f(h)=h^{\{1\}}\sstac{T}
h^{\{2\}}$ for all $h\in H$, and compute
\begin{eqnarray*}
{\widetilde {\rm can}}_T\big(\ell_T(h)\big)
&=& \conv{j}(h\di){h\dii}^{\{1\}}{h\dii}^{\{2\}}j(h_{(3)})\tz
\stac{R} j(h_{(3)})\ti\\
&=& \conv{j}(h\di) \eta_L\big( \pi_L(h\dii)\big)j({h_{(3)}}\ui)
\stac{R}{h_{(3)}}\uii\\
&=& \conv{j}(h\di)j({h_{(2)}}\ui)\stac{R}{h_{(2)}}\uii
= \conv{j}({h\ui}\di)j({h\ui}\dii)\stac{R}{h}\uii\\
&=& \eta_R\big(\pi_R(h\ui)\big)\stac{R} h\uii
=1_A\stac{R} h,
\end{eqnarray*}
where the first equality follows by \eqref{eq:tcan} and the fact
that the range of $f$ is in $B\sstac{T} B$. The second equality
follows by the hypothesis $\mu_B\circ f=\eta_L\circ \pi_L$ and the
right $\hH$-colinearity of $j$. The third one follows by the left
$L$-linearity of $j$ (i.e.\ \eqref{eq:js}) and the counit property
of $\pi_L$. The fourth equality follows by the left
$\hH_L$-colinearity of $\gamma_R$. In the penultimate step we used
that $\conv{j}$ is a left convolution inverse of $j$, i.e.\
\eqref{eq:linv}. This means that $\ell_T$ given in
\eqref{eq:strcon} is a strong $T$-connection and completes the
proof of the theorem.
\end{proof}

Take a bijective entwining structure $({A},{C},\psi)_R$ over a
non-commutative base algebra $R$ and consider a $T$-flat entwined
extension $B\subseteq A$ in the sense of \cite[Definition
5.2]{BohBrz:str} ($T$ is a subalgebra of $B$). Given a strong
$T$-connection in $B\subseteq A$, one constructs
 a family of maps of Abelian groups from the
Grothendieck group of ${C}$-comodules to the even $T$-relative cyclic
homology groups of $B$ (cf.\  \cite[Theorem 5.4]{BohBrz:str}).
 This family of maps is termed the {\em T-relative Chern-Galois
character}.
Comodule algebras for Hopf algebroids (with bijective
antipodes) provide
examples of (bijective) entwining structures over non-commutative
bases, hence the
general theory worked out in \cite{BohBrz:str} can be applied to such
algebras. In
particular, the components of the $T$-relative Chern-Galois characters,
corresponding to
strong $T$-connections in Theorem~\ref{thm.class.con} for a $T$-flat cleft
extension of a Hopf algebroid with bijective antipode, have been computed in
\cite[Example 5.6]{BohBrz:str}.

It is important to note, however, that the $T$-relative Chern-Galois
character, a priori, depends on the choice of a strong $T$-connection. Its
independence is proven in \cite[Theorem 5.14]{BohBrz:str}, under the
assumption that the $T$-flat entwined extension $B\subseteq A$
enjoys the following properties:
\begin{blist}
\item $A$ is a locally projective right $T$-module;
\item the extension $B\subseteq A$ splits as a $B$-$T$ bimodule.
\end{blist}
In the remainder of this section we analyse the meaning of these conditions
and
of the $T$-flatness in the case
of cleft Hopf algebroid extensions. In this way we find sufficient conditions
for
the existence and the strong-connection-independence of the relative
Chern-Galois character computed in \cite[Example 5.6]{BohBrz:str}.
\begin{definition}\label{def.tot.int}
Let  $B\subseteq A$ be a right extension of a Hopf algebroid $\hH$ with a
bijective antipode and let $T$ be a subalgebra of $B$. View $A$ as a
left $\hH_R$-comodule with coaction \eqref{eq:Aro}
{\color{blue}
and $\hH_L$-coaction $a\mapsto S^{-1}(a^{[1]})\ot_L a^{[0]}$}.
A {\em left total $T$-integral} is a left
{\color{blue} 
$\hH$}-colinear map $\vartheta: H\to A$ such that $\vartheta(H)\subseteq
A^T:= \{\ a\in A \; |\; ta=at\quad \forall t\in T\ \}$
and $\vartheta(1_H) = 1_A$.
\end{definition}
For example, the convolution inverse of a normalised cleaving map is a left
total $k$-integral by Lemma~\ref{lem:roaj}.
By arguments similar  to those used to prove Lemma~\ref{lem:coinv}, any left
total $T$-integral $\vartheta$ determines a $B$-$T$ bilinear section of the
extension $B\subseteq A$,
\begin{equation}\label{eq:intsplit}
a\mapsto a\tz\vartheta(a\ti).
\end{equation}
The next proposition shows that, for a cleft extension of a Hopf
algebroid with a bijective antipode,
this is a one-to-one  correspondence.

{\color{blue}
Note that if $B\subseteq A$ is a cleft extension by a Hopf algebroid $\hH$
with a bijective antipode, then 
similarly to the proof of Lemma \ref{lem:H-H_R}, any morphism $f\in
\Hom^{\hH_R,-}(H,A)$ satisfies $f(h)=\conv{j}(h_{(1)}) 
j \big(S^{-1}(f(h_{(2)})_{[1]}) \big)f(h_{(2)})_{[0]}$, for any $h\in H$,
where $j$ is a cleaving map with convolution inverse $\conv{j}$. 
Hence $f$ is also left $\hH_L$-colinear,
with respect to the $\hH_L$-coaction $a\mapsto S^{-1}(a^{[1]})\ot_L a^{[0]}$
on $A$. In particular, in this situation a left total $T$-integral is the same
as a left $\hH_R$-colinear map $\vartheta: H\to A$ such that
$\vartheta(H)\subseteq A^T$ and $\vartheta(1_H) = 1_A$.}

\begin{proposition}\label{prop.tot.int}
Let  $B\subseteq A$ be a cleft extension of a Hopf algebroid $\hH$ with a
bijective antipode, and let $T$ be a
subalgebra of $B$. Then $B$-$T$ bilinear sections of the extension $B\subseteq
A$ are in bijective correspondence with left total $T$-integrals in
$B\subseteq A$.
\end{proposition}
\begin{proof}
Let $j$ be a cleaving map for $B\subseteq A$. In terms of $j$ we
construct the inverse of the map associating to a left total
$T$-integral $\vartheta$ the section \eqref{eq:intsplit}. To a
$B$-$T$ bilinear section $\varphi$, associate the map
\begin{equation}\label{eq:splitint}
\vartheta:H\to A,\qquad h\mapsto \conv{j}(h\di)\varphi(j(h\dii)).
\end{equation}
Since $j(1_H)$ is an element of $B$, $\vartheta(1_H)=
\conv{j}(1_H)j(1_H)\varphi(1_A)=1_A$. The left 
{\color{blue}
$\hH$}-colinearity
of $\vartheta$ follows by the left 
{\color{blue}
$\hH$}-colinearity of
$\gamma_L$ and $\conv{j}$, and the fact that the range of
$\varphi$ is equal to $B$. It remains to check that the range of
$\vartheta$ is in $A^T$. Note that by its left $B$-linearity,
$\varphi$ is determined by the left $L$-linear map $\varphi\circ
j:H\to B$. Indeed, for all $a\in A$,
$$
\varphi(a) =\varphi\big( a\tz \conv{j}({a\ti}\di)
j({a\ti}\dii)\big) = \varphi\big( {a\bz}\tz \conv{j}({a\bz}\ti)
j({a\bi})\big) ={a\bz}\tz \conv{j}({a\bz}\ti)\varphi(j(a\bi)),
$$
where in the last equality Lemma \ref{lem:coinv} has been used.
Hence the right $T$-linearity of $\varphi$ is equivalent to
\begin{equation}\label{eq:Tlin}
a\tz t \conv{j}({a\ti}\di)\varphi(j({a\ti}\dii))=\varphi(a)
t,\qquad \textrm{for all } a\in A \textrm{ and } t\in T.
\end{equation}
Take any $h\in H$ and apply \eqref{eq:Tlin} to $a=j(h)$. By the right
$\hH_R$-colinearity of $j$,
$$
j(h\ui) t \conv{j}({h\uii}\di) \varphi(j({h\uii}\dii ))= \varphi(j(h))
t.
$$
Hence
$$
\vartheta(h) t =\conv{j}(h\di) j({h\dii}\ui)t
\conv{j}({{h\dii}\uii}\di)\varphi(j({{h\dii}\uii}\dii))
=\eta_R(\pi_R(h\ui)) t \vartheta(h\uii)=t \vartheta(h),
$$
where the second equality follows by the Hopf algebroid axiom
\eqref{ax:cpcolin} and the last one follows by the fact that the
elements of $B$ (and hence, in particular, the elements of $T$)
commute with $\eta_R(r)$, for $r\in R$, and by the left
$R$-linearity of $\vartheta$. It is checked by a routine
computation that the map, associating to a $B$-$T$ bilinear
section $\varphi$ of the inclusion $B\subseteq A$ the left total
$T$-integral \eqref{eq:splitint}, is the inverse of the map,
associating to the left total $T$-integral $\vartheta$ the $B$-$T$
bilinear section \eqref{eq:intsplit}.
\end{proof}

For a cleft extension $B\subseteq A$ of a Hopf algebroid $\hH$ with a
bijective antipode, consider the $B$-$B$ bilinear map,
$$
A\to A\stac{R} H,\qquad a\mapsto a\tz\stac{R} a\ti-a\stac{R} 1_H,
$$
where $A\sstac{R} H$ is a $B$-$B$ bimodule via the first
tensorand. For any subalgebra $T$ of $B$ this map projects to the
map $\upsilon_T: A/[A,T]\to A/[A,T]\sstac{R} H$, where $[A,T]=\{\
\sum_k (a_k t_k-t_ka_k)\ |\ a_k\in A,\ t_k\in T\ \}$ is a right
$R$-submodule of $A$. Following \cite[Definition 5.2]{BohBrz:str},
the extension $B\subseteq A$ is said to be {\em $T$-flat} if $B$
and $A$ are flat left and right $T$-modules and the obvious map
\begin{equation}\label{eq:(5.1)}
B/[B,T]\to {\rm ker}\ \upsilon_T, \qquad [b]_B\mapsto [b]_A,
\end{equation}
(where $[\ ]_B$ denotes the equivalence class in $B/[B,T]$ and $[\ ]_A$
denotes the equivalence class in $A/[A,T]$) is an isomorphism.
\begin{proposition}\label{prop.T-flat}
Let $\hH$ be a Hopf algebroid with a bijective antipode.
A cleft $\hH$-extension  $B\subseteq A$
which splits as a $B$-$T$ bimodule for some subalgebra $T$
of $B$, is $T$-flat if and only if $A$ is a flat left and right
$T$-module.
\end{proposition}
\begin{proof}
Since a direct summand of a flat module is flat, it suffices to
prove that the existence of a $B$-$T$ bimodule splitting of the
inclusion, i.e.\ the existence of a left total $T$-integral,
implies that the map \eqref{eq:(5.1)} is an isomorphism. In order
to prove injectivity of \eqref{eq:(5.1)}, choose $b\in B$ such
that $[b]_A=0$. This means the existence of finite sets $\{a_k\}$
in $A$ and $\{ t_k\}$ in $T$ such that $b=\sum_k (a_k t_k-t_k
a_k)$. Applying a $B$-$T$ bilinear section $\varphi$ of the
extension $B\subseteq A$ to this identity, we obtain $b=\sum_k
(\varphi(a_k) t_k- t_k \varphi(a_k))$, hence $[b]_B=0$. In order
to prove the surjectivity of the map \eqref{eq:(5.1)}, choose
$a\in A$ such that $\upsilon_T([a]_A)=0$. This means the existence
of finite sets $\{a_k\}$ in $A$, $\{h_k\}$ in $H$ and $\{ t_k\}$
in $T$ such that
\begin{equation}\label{eq:kerup}
a\tz\stac{R} a\ti-a\stac{R} 1_H=\sum_k (a_k t_k\stac{R} h_k -t_k a_k\stac{R}
h_k).
\end{equation}
By Proposition~\ref{prop.tot.int}, there is a left total
$T$-integral $\vartheta$ in $B\subseteq A$. Apply $\mu_A\circ(A\sstac{R}
\vartheta)$ to \eqref{eq:kerup} to obtain
$$
{a}\tz \vartheta({a\ti})- a=
\sum_k \left( a_k \vartheta(h_{k}) t_k-
       t_k a_k \vartheta(h_{k})\right).
$$
This proves that $[a]_A=[{a}\tz\vartheta({a\ti})]_A$. Since
${a}\tz\vartheta({a\ti})$ is an element of $B$,
$[a]_A$ belongs to the image of the map \eqref{eq:(5.1)}.
\end{proof}

Combining
Proposition~\ref{prop.tot.int} and Proposition~\ref{prop.T-flat} with
\cite[Theorem~5.14]{BohBrz:str} we obtain
\begin{corollary}\label{cor:cleftChG}
Let $B\subseteq A$ be a cleft extension of a Hopf algebroid with a bijective
antipode and let $T$ be a subalgebra of $B$.
Assume that:
\begin{blist}
\item $A$ is a flat left $T$-module and a locally projective right $T$-module;
\item there exists a left total $T$-integral
for the extension $B\subseteq A$;
\item there exists a strong $T$-connection.
\end{blist}
Then there exists a $T$-relative Chern-Galois character, independent on the
choice of the strong $T$-connection in (c).
\end{corollary}
We close the section with some examples of Hopf algebroid cleft extensions,
in which
there exist (strong-connection-independent) relative Chern-Galois
characters.
\begin{example}
Let $B\subseteq A$ be a cleft extension of a Hopf algebroid with a
bijective antipode and $T$ a {\em separable} $k$-subalgebra of $B$.
In light of \cite[Proposition 3.2 (1)]{BohBrz:str}, since $B\subseteq A$
splits as a left $B$-module by Lemma \ref{lem:dirsum}, it splits as a
$B$-$T$ bimodule. The corresponding total $T$-integral $\vartheta$
in $B\subseteq A$ is
$$
\vartheta(h) =  \sum_ie_i\conv{j}(h)j(1_H)e^i,
$$
where $\sum_i e_i\sstac{k}e^i \in T\sstac{k}T$ is a separability
idempotent. Therefore, if $A$ is a flat left $T$-module and a
locally projective right $T$-module and there exists a strong
$T$-connection $\ell_T$, then there exists a corresponding
$T$-relative Chern-Galois character which is independent of
$\ell_T$.
\end{example}
\begin{example}
The base algebra $R$ of a Hopf algebroid $\hH$ is a right
$\hH$-comodule algebra with 
{\color{blue}
$\hH_R$-coaction 
$r\mapsto 1_R \ot_R s_R(r)$
and $\hH_L$-coaction 
$r\mapsto 1_R \ot_L s_R(r)$}. 
It follows by
\cite[Theorem 3.2]{Bohm:hgdint} that a Hopf algebroid $\hH$ with a
bijective antipode is coseparable (as an $L$- or, equivalently, as
an $R$-coring) if and only if there exists a left total
$k$-integral $\lambda$ for the $\hH$-extension $I\subseteq R$,
where $I$ is the 
{\color{blue}
$\hH_R$}-coinvariant subalgebra of $R$, $I=\{\ r\in
R\ |\ s_R(r)=t_R(r)\ \}$. Note that if $\hH$ is a coseparable Hopf
algebroid, then any right $\hH$-extension $B\subseteq A$ is split
as a $B$-$B$ bimodule by the map
$$
A\to B,\qquad a\mapsto a\tz \eta_R(\lambda(a\ti)).
$$
Let $B\subseteq A$ be a cleft extension of a {\em coseparable} Hopf
algebroid with
a bijective antipode and $T$ a subalgebra of $B$. Then there is a
left total $T$-integral $\vartheta = \eta_R\circ \lambda$. If $A$ is a locally
projective right $T$-module and a flat left $T$-module and there exists a
strong $T$-connection $\ell_T$,
then the corresponding $T$-relative Chern-Galois character  is
independent of $\ell_T$.
\end{example}

\appendix
\section{Weak cleft extensions and weak crossed products}

Many of the results described in Sections~\ref{sec:Hcleft},
\ref{sec:crossp} and \ref{sec.chern}
(e.g. Lemma~\ref{lem:dirsum}, Corollary~\ref{cor:ff},
Corollary~\ref{cor:rel.inj} or Theorem~\ref{thm.str.con.can},
Corollary~\ref{cor:cleftChG})
remain valid if the right $\hH$-extension
$B\subseteq A$ is an $\hH_R$-Galois
extension but, instead of the normal basis condition
Theorem~\ref{thm.main.cleft}~(2)(b),
 $A$ is only a direct summand of $B\sstac{L} H$ as
a left $B$-module and as a right $\hH$-comodule. Such extensions can be
studied along the same lines as in
Sections \ref{sec:Hcleft}, \ref{sec:crossp} and \ref{sec.chern}.
In this appendix we present the results of such studies; we give no proofs as
these are very similar to the proofs of corresponding results in
preceding sections.

Motivated by the forthcoming analogue of Theorem~\ref{thm.main.cleft}
(Theorem~\ref{thm:wcleft}), we introduce the following
weakening of Definition~\ref{def:cleft}.
\begin{definition}\label{def:wcleft}
Let $\hH$ be a Hopf algebroid. A right $\hH$-extension $B\subseteq A$ is
{\em weak cleft} if
\begin{blist}
\item in addition to its canonical $R$-ring structure, $A$ possesses
an $L$-ring structure and $B$ is an $L$-subring of $A$;
\item there exists a left $L$-linear right $\hH$-colinear morphism $j:H\to A$,
  with left convolution inverse $\lconv{j}$, which is right
{\color{blue}
$\hH$}-colinear in the sense of identities
  \eqref{eq:third} and \eqref{eq:roaj}.
\end{blist}
A map $j$, satisfying condition (b), is called a {\em weak cleaving map}.
\end{definition}
Note that in the situation described in Definition~\ref{def:wcleft},
the assumption that $\lconv{j}$ satisfies \eqref{eq:roaj} implies that
the image of the map $A\to A$, $a\mapsto
a\tz \lconv{j} (a\ti )$ is contained in $B$.
Hence a weak $\hH$-cleft extension $B\subseteq A$ is split by the left
$B$-linear map
\eqref{eq:split} after replacing $\conv{j}$ with $\lconv{j}$.
\begin{theorem}\label{thm:wcleft}
Let $\hH$ be a Hopf algebroid and $B\subseteq A$ a right $\hH$-extension.
Then the following statements are equivalent:
\begin{zlist}
\item $B\subseteq A$ is a weak $\hH$-cleft extension.
\item
\begin{blist}
\item The extension $B\subseteq A$ is $\hH_R$-Galois; \item $A$ is
a direct summand of $B\sstac{L} H$ as a left $B$-module and
  right $\hH$-comodule.
\end{blist}
\end{zlist}
\end{theorem}
In particular, Theorem~\ref{thm:wcleft} implies that if
$\hH$ is projective as a left $L$-module then, for any weak cleft
$\hH$-extension $B\subseteq A$, $A$ is a faithfully flat left $B$-module.

Recall that in Section~\ref{sec:crossp} we  applied a
(unnormalised) gauge transformation to a general
cleaving map in order to normalise it as $j(1_H)=1_B=\conv{j}(1_H)$. However,
in
the case when $j$ possesses a left convolution inverse $\lconv{j}$ only,
there is no guarantee for $j(1_H)$ to be an invertible element of $B$. Hence
it can not be gauge transformed to the unit element in general. The need to
describe this more general situation leads to the following generalisations of
Definitions \ref{def:meas} and \ref{def:cocycle}.
\begin{definition}\label{def:wmeas}
Let ${\mathcal L}=(H,L,s,t,\gamma,\pi)$ be a left bialgebroid and
$B$ an $L$-ring. ${\mathcal L}$ {\em weakly measures} $B$ if there
exits a $k$-linear map, termed a {\em weak measuring},
$H\sstac{k} B\to B$,  $h\sstac{k} b\mapsto h\cdot b$ that
satisfies properties (b) and (c) in Definition \ref{def:meas}.
\end{definition}
\begin{definition}\label{def:wcoc}
Let ${\mathcal L}=(H,L,s,t,\gamma,\pi)$ be a left bialgebroid and
$B$ an $L$-ring, weakly measured by ${\mathcal L}$. A $B$-valued
{\em weak 2-cocycle} $\sigma$ on ${\mathcal L}$ is a $k$-linear
map $H\sstac{L^{op}} H\to B$ (where the right and left
$L^{op}$-module structures on $H$ are given by right and left
multiplication by $t(l)$, respectively) satisfying properties (a),
(b) and (d) in Definition \ref{def:cocycle} and, in addition, for
all $h,k\in H$,
\begin{equation}\label{eq:wcoc}
\sigma(h\di,k\di)(h\dii k\dii \cdot 1_B)=\sigma(h,k).
\end{equation}
A weakly ${\mathcal L}$-measured $L$-ring $B$ is called a
{\em $\sigma$-twisted left ${\mathcal L}$-module} if a weak $2$-cocycle
$\sigma$ satisfies property (f) in Definition~\ref{def:cocycle} and there
exist elements $x$ and ${\tilde x}$ in $B$ such that, for all $b\in B$ and
$h\in H$,
\begin{eqnarray*}
{\tilde x} x=1_B\qquad &\textrm{and}&\qquad x b {\tilde x}= 1_H\cdot b,\\
\sigma(1_H,h)=x(h\cdot 1_B)\qquad &\textrm{and}&\qquad \sigma(h,1_H)=h\cdot x.
\end{eqnarray*}
\end{definition}
It is easy to see that a $B$-valued 2-cocycle is also a weak 2-cocycle. If the
$L$-ring $B$ is a $\sigma$-twisted $\cL$-module for a 2-cocycle $\sigma$, then
it is a $\sigma$-twisted $\cL$-module also in the weaker sense of
Definition~\ref{def:wcoc} with $x=1_B={\tilde x}$.

Recall from \cite[p.\ 39]{CaeDeG:mod} that, for a non-unital ring $A$,
an element
$e\in A$ such that, for all $a\in A$, $ea=ae=ae^2$ is called a {\em
preunit}. Proposition~\ref{lem:excp} can be extended to the case of
weak cocycles as follows.
\begin{proposition}\label{prop:wcp}
Let ${\mathcal L}=(H,L,s,t,\gamma,\pi)$ be a left bialgebroid and
$B$ an $L$-ring, weakly measured by ${\mathcal L}$. Let $\sigma:
H\sstac{L^{op}} H\to B$ be a map that satisfies properties (a) and
(b) in Definition~\ref{def:cocycle} and condition \eqref{eq:wcoc}
in Definition~\ref{def:wcoc}. Consider the $k$-module $B\sstac{L}
H$ in Proposition~\ref{lem:excp}, and the following assertions.
\begin{blist}
\item $B\sstac{L} H$ is an associative (possibly non-unital)
algebra with
  multiplication \eqref{sigma.prod}.
\item There exists ${\tilde y}\in B$ such that ${\tilde y}
\sigma(1_H,h) = h \cdot 1_B$, for all $h\in H$, and ${\tilde y}
\sstac{L} 1_H$ is a preunit for the algebra in part (a) (hence
$A:=\{(b\sstac{L} h)({\tilde y}\sstac{L} 1_H)\; |\; b\sstac{L}
h\in B\sstac{L} H\}$ is a right ${\mathcal L}$-comodule algebra
with coinvariant subalgebra $\{ (b \sstac{L} 1_H)({\tilde y}
\sstac{L} 1_H)\; |\; b\in B \}$). \item The map $B\to A^{co\cL}$,
$b\mapsto (b\tilde{y}\sstac{L}1_H)(\tilde{y}
  \sstac{L} 1_H)$ is an algebra isomorphism.
\end{blist}
These assertions hold if and only if $\sigma$ is a weak 2-cocycle and $B$ is a
$\sigma$-twisted left $\cL$-module. In this case $A$ is called a {\em weak
  crossed product} of $B$ with $\cL$ with respect to the weak 2-cocycle
$\sigma$.
\end{proposition}
Our next task is to characterise equivalent weak crossed products, in analogy
with Theorem~\ref{thm.equiv}.
\begin{definition}\label{def:wcpeq}
Let ${\mathcal L}=(H,L,s,t,\gamma,\pi)$ be a left bialgebroid and
$B$ an $L$-ring. Weak crossed product algebras of $B$ with $\cL$ are said to
be {\em equivalent} if they are isomorphic via a left $B$-linear isomorphism
of right $\cL$-comodule algebras.
\end{definition}
Note that a (left $B$-linear right $\cL$-colinear) isomorphism of
weak crossed product algebras of $B$ with $\cL$ in
Definition~\ref{def:wcpeq} needs not extend to the (non-unital)
algebra $B\sstac{L} H$. The following lemma extends
Corollary~\ref{cor.gauge}.
\begin{lemma}\label{lem:wgauge}
Let ${\mathcal L}=(H,L,s,t,\gamma,\pi)$ be a left bialgebroid and
$B$ an $L$-ring, weakly measured by ${\mathcal L}$. Let $\sigma$ be
a weak 2-cocycle, such that $B$ is a $\sigma$-twisted left $\cL$-module. Let
$\chi$ and ${\tilde \chi}$ be morphisms in ${\rm Hom}_{L,L}(H,B)$ such that,
for all $h\in H$,
\begin{eqnarray}
&&{\tilde \chi}(h\di)\chi(h\dii)=h\cdot 1_B\qquad \textrm{and}\label{eq:wga}\\
&&{\tilde \chi}(h\di)\chi(h\dii){\tilde \chi}(h_{(3)})={\tilde \chi}(h),
\quad {\chi}(h\di){\tilde \chi}(h\dii){\chi}(h_{(3)})={\chi}(h).\label{eq:wgc}
\end{eqnarray}
Then \eqref{gauge.mea} defines a weak measuring
of $\cL$ on $B$ and \eqref{gauge.coc} defines a weak 2-cocycle
$\sigma^\chi$, such that $B$ is a $\sigma^\chi$-twisted left $\cL$-module.
\end{lemma}
A pair $\chi, {\tilde \chi}\in {\rm Hom}_{L,L}(H,B)$, satisfying
(\ref{eq:wga}- \ref{eq:wgc}), is called a {\em gauge transformation} of the
weak crossed product of $B$ with ${\mathcal L}$.
Gauge transformations form a groupoid,
with multiplication, the convolution product $\diamond$ in the first
component, and its opposite in the second one.
The left unit of a gauge transformation $(\chi, {\tilde \chi})$ is
$({\chi}\diamond {\tilde \chi}, \chi\diamond {\tilde \chi})$
and its right unit is
$({\tilde \chi}\diamond {\chi}, {\tilde \chi}\diamond {\chi})$.
The inverse of $(\chi,{\tilde \chi})$ is $({\tilde \chi}, {\chi})$.
\begin{theorem}\label{thm:wcpeq}
Let ${\mathcal L}=(H,L,s,t,\gamma,\pi)$ be a left bialgebroid and
$B$ an $L$-ring. Two weak crossed products of $B$ with $\cL$ are
equivalent if
and only if they are related by a gauge transformation.
\end{theorem}
In order to make connection between weak crossed products and weak cleft
extensions, the notion of invertible weak 2-cocycles is needed.
\begin{definition}\label{def:wcocinv}
Let ${\mathcal L}=(H,L,s,t,\gamma,\pi)$ be a left bialgebroid and
$B$ an $L$-ring weakly measured by ${\mathcal L}$. An {\em
inverse} of a $B$-valued weak 2-cocycle $\sigma$ on $\cL$ is a
$k$-linear map, ${\tilde \sigma}:H\sstac{L} H\to B$ (where the
right and left $L$-module structures on $H$ are given by right and
left multiplication by $s(l)$, respectively) satisfying properties
(a), (b) and (c) in Definition~\ref{def:cocinv} and, in addition,
for all $h,k\in H$,
$$
(h\di k\di\cdot 1_B){\tilde \sigma}(h\dii, k\dii)={\tilde \sigma}(h,k).
$$
\end{definition}
If $\sigma$ is a 2-cocycle in the sense of Definition
\ref{def:cocycle} (in particular the measuring satisfies also property
(a) in Definition~\ref{def:meas}), then
Definition~\ref{def:wcocinv} is equivalent to Definition \ref{def:cocinv}.
By an argument similar to the proof of Lemma~\ref{lem:tsinorm}, the
convolution inverse of a weak 2-cocycle is unique, provided it exists.
A generalisation of Theorems~\ref{thm.cross.cleft} and \ref{thm:crossed}
is given in the following
\begin{theorem}\label{thm:wclcp}
Let $\hH$ be a Hopf algebroid. A right $\hH$-extension  $B\subseteq A$
is weak $\hH$-cleft if and only if $A$ is isomorphic to a weak crossed
product of $B$ with the constituent left bialgebroid $\hH_L$ of $\hH$, with
respect to an invertible weak 2-cocycle.
\end{theorem}
Analogously to Corollary~\ref{cor.gauge.cleft}, Theorems~\ref{thm:wclcp} and
\ref{thm:wcpeq} lead to the following
\begin{corollary}
Let $\hH$ be a Hopf algebroid and $B\subseteq A$ a weak
$\hH$-cleft extension. Let $j:H\to A$ be  a weak cleaving map with
left convolution inverse $\lconv{j}$, satisfying conditions
\eqref{eq:third} and \eqref{eq:roaj}. Then the map $h\cdot b\colon
= j(h\ui)b \lconv{j}(h\uii)$, for all $b\in B$ and $h\in H$, is a
weak measuring of the constituent left bialgebroid $\hH_L$ on $B$.
A map $j':H\to A$ is a weak cleaving map if and only if there
exist morphisms $\chi, {\tilde \chi}\in {\rm Hom}_{L,L}(H,B)$,
satisfying (\ref{eq:wga}-\ref{eq:wgc}), in terms of which
$j':h\mapsto \chi(h\di) j(h\dii)$.
\end{corollary}

Let $\hH$ be a Hopf algebroid with a bijective antipode and let
$B\subseteq A$ be a weak $\hH$-cleft extension and $T$ a
$k$-subalgebra of $B$. Let $j$ be a weak cleaving map with a left
convolution inverse $\lconv{j}$, satisfying conditions
\eqref{eq:third} and \eqref{eq:roaj}. Any morphism $f\in \lrhom
LLH{B\sstac{T} B}$ such that, for all $h\in H$,
$\mu_B\big(f(h)\big)=j(h\ui)\lconv{j}(h\uii)$, determines a strong
$T$-connection via \eqref{eq:strcon}. Conversely, any strong
$T$-connection is of this form (though the correspondence
\eqref{eq:strcon} is not bijective in the weak case).

Any $B$-$T$ bimodule section of a weak cleft Hopf algebroid extension
$B\subseteq A$, for a subalgebra $T$ of $B$,
corresponds to a left total $T$-integral via \eqref{eq:intsplit}
(although the correspondence
between $B$-$T$ sections and total integrals
is not bijective in the weak case). 
Hence Corollary \ref{cor:cleftChG}
is valid without modification for weak cleft extensions of Hopf
algebroids with bijective antipode.
\bigskip

A weak Hopf algebra $(W,\Delta,\epsilon,S)$ determines a 
Hopf algebroid ${\mathcal W}$ with constituent left bialgebroid
${\mathcal W}_L$ over the `left' subalgebra $W^L$ of $W$, right
bialgebroid ${\mathcal W}_R$ over the `right' subalgebra $W^R$,
and antipode $S$. The category of right comodules for the
coalgebra $(W,\Delta,\epsilon)$ is isomorphic to the category of
right ${\mathcal W}$-comodules as a monoidal category. As a
consequence, also the respective notions of comodule algebras and
of coinvariants are equivalent (cf. \cite{BrCaMi:DHalg}). Let $A$
be a right $W$- (or, equivalently, ${\mathcal W}$-) comodule
algebra with coinvariants $B$. By \cite[Theorem
2.11]{AlAl:wcleftgal}, the extension $B\subseteq A$ is $W$-cleft
(i.e. the corresponding weak entwining structure is cleft in the
sense of \cite[Definition 1.9]{AlAl:wcleft}) if and only if it is
$W$-Galois and $A$ is a direct summand of $B\sstac{k} W$ as a left
$B$ module right $W$-comodule. By \cite[Example 3.5]{Bohm:gal} the
$W$-Galois property is equivalent to the ${\mathcal W}_R$-Galois
property, hence Theorem~\ref{thm:wcleft} implies that any weak
${\mathcal W}$-cleft extension is weak $W$-cleft (but not
conversely).

\section*{Acknowledgements}
GB is grateful to Korn\'el Szlach\'anyi and Joost Vercruysse for
inspiring discussions. Her work is supported by the Hungarian
Scientific Research Fund OTKA T 043 159 and the Bolyai J\'anos
Fellowship.

\end{document}